# Analytical Methods for Squaring the Disc

Chamberlain Fong
spectralfft@yahoo.com
Seoul ICM 2014

**Abstract** – We present and discuss several old and new methods for mapping a circular disc to a square. In particular, we present analytical expressions for mapping each point (u,v) inside the circular disc to a point (x,y) inside a square region. Ideally, we want the mapping to be smooth and invertible. In addition, we put emphasis on mappings with desirable properties. These include conformal, equiareal, and radially-constrained mappings. Finally, we present applications to logo design, panoramic photography, and hyperbolic art.

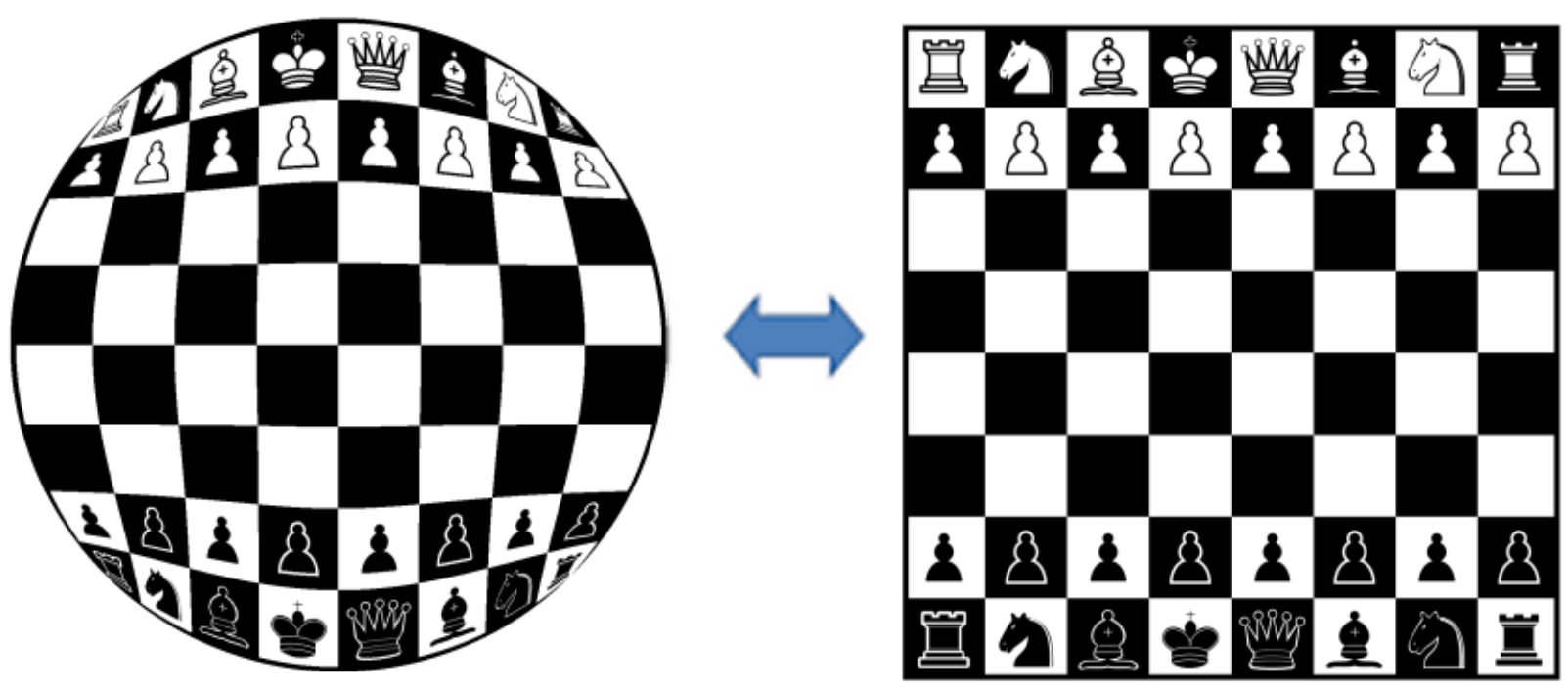

## 1 Introduction

The circle and the square are among the most common shapes used by mankind. It is certainly worthwhile to study the mathematical correspondence between the two. In this paper, we shall discuss ways to map a circular region to a square region and back. There are infinitely many ways of doing this mapping. Of particularly interest to us are mappings with nice closed-form invertible equations. We emphasize the importance of invertible equations because we want to perform the mapping back and forth between the circular disc and the square. We shall present and discuss several such mappings in this paper.

## 1.1 Organization of this Paper

We anticipate that there will be two types of people who might read this paper. The 1st kind would be those who just want to get the equations to map a circular region to a square; and do not really care about proofs or derivations. The 2nd kind would be more interested in the mathematical details behind the mappings. Therefore, we shall organize this paper into two parts. The 1st part will only contain equations for mapping a circular disc to a square region and back. The 2nd part will delve more into mathematical details and discuss some desirable properties of the different mappings. In addition to this, we will also discuss real-world applications of these mappings.

There are accompanying presentation slides and C++ implementation of this paper available on these websites:
http://www.slideshare.net/chamb3rlain/mappings-for-squaring-the-circular-disc
http://squircular.blogspot.com

## 1.3 Canonical Mapping Space

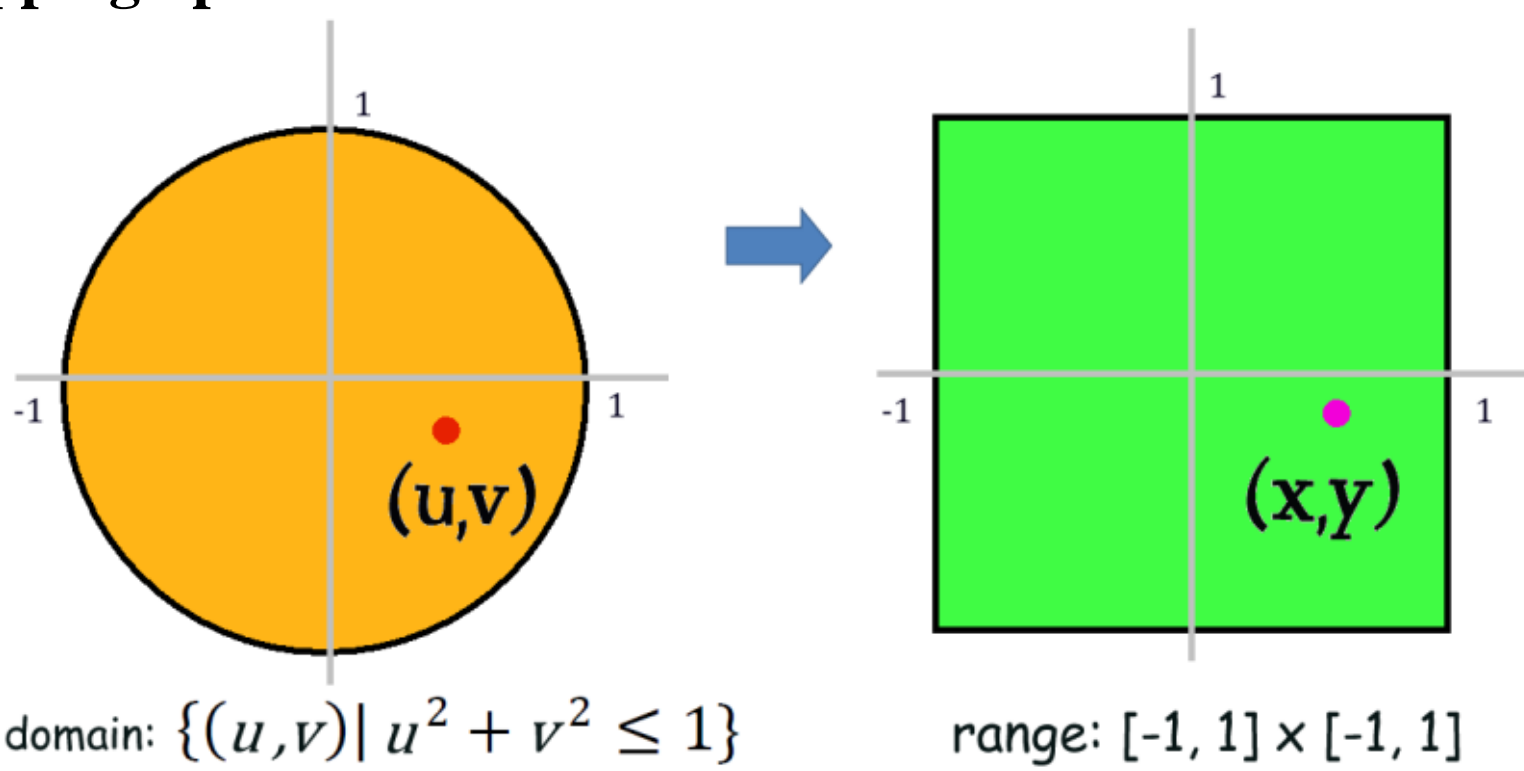

The canonical space for the mappings presented here is the unit disc centered at the origin with a square circumscribing it. This unit disc is defined as the set $\mathcal{D} = \{(u,v)|\ u^2 + v^2 \leq 1\}$. The square is defined as the set $\mathcal{S} = [-1,1] \text{ x } [-1,1]$. This square has a side of length 2. We shall denote (u,v) as a point in the interior of the unit disc and (x,y) as the corresponding point in the interior of the square after the mapping. In this paper, we shall present several equations that relate (u,v) to (x,y).

Mathematically speaking, we want to find functions *f* that maps every point (u,v) in the circular disc to a point (x,y) in the square region and vice versa. In others word, we want to derive equations for *f* such that (u,v) =*f*(x,y) and (x, y) = $f^{\,1}$(u, v) .

In addition, we shall impose these three conditions on the mapping:

- (0,0) = *f*(0,0) , i.e. the center of the circle corresponds to the center of the square
- (±1,0) = *f*(±1,0), i.e. the extreme points of the shapes in the x-axis match
- (0, ±1) = *f*(0, ±1), i.e. the extreme points of the shapes in the y-axis match

We shall denote the 1st condition as the *canonical central constraint* of the mapping. Also, we shall denote the 2nd and 3rd conditions as the *canonical axial constraints* of the mapping.

## 1.4 A Classic Problem in a Modern Guise?

The mapping of the circular disc to a square region is similar but not equivalent to the classic mathematical problem of "squaring the circle". For one thing, in the classic mathematical problem, one is restricted to using only a straightedge and a compass. Our problem concerns finding mapping equations that a computer can calculate. In particular, we want explicit equations to mapping each point *(u,v)* to point *(x,y)*. The two problems are superficially similar but ultimately quite different. One problem has to do with geometric construction; while the other problem has to do with finding a two-dimensional mapping function.

## 1.5 The Mappings

In the next four pages, we shall present four mappings for converting the circular disc to a square and vice versa. We include pictures of a radial grid inside the circle converted to a square; and a square grid converted to a disc. This is followed by equations for the mappings. In these equations, we shall make use a common math function called the signum function, denoted as *sgn(x).* The signum function is defined as

$$\mathrm{sgn(x)} = \frac{|\mathrm{x}|}{\mathrm{x}} = \begin{cases} -1 & \text{if x} < 0, \\ 0 & \text{if x} = 0, \\ 1 & \text{if x} > 0. \end{cases}$$

Also note that for the sake of brevity, we have not singled out cases when there are divisions by zero in the mapping equations. For these special cases, just equate *x=u*, *y=v* and vice versa when there is an unwanted division by zero in the equations. This usually happens when *u=0* or *v=0* or both.

# Simple Stretching

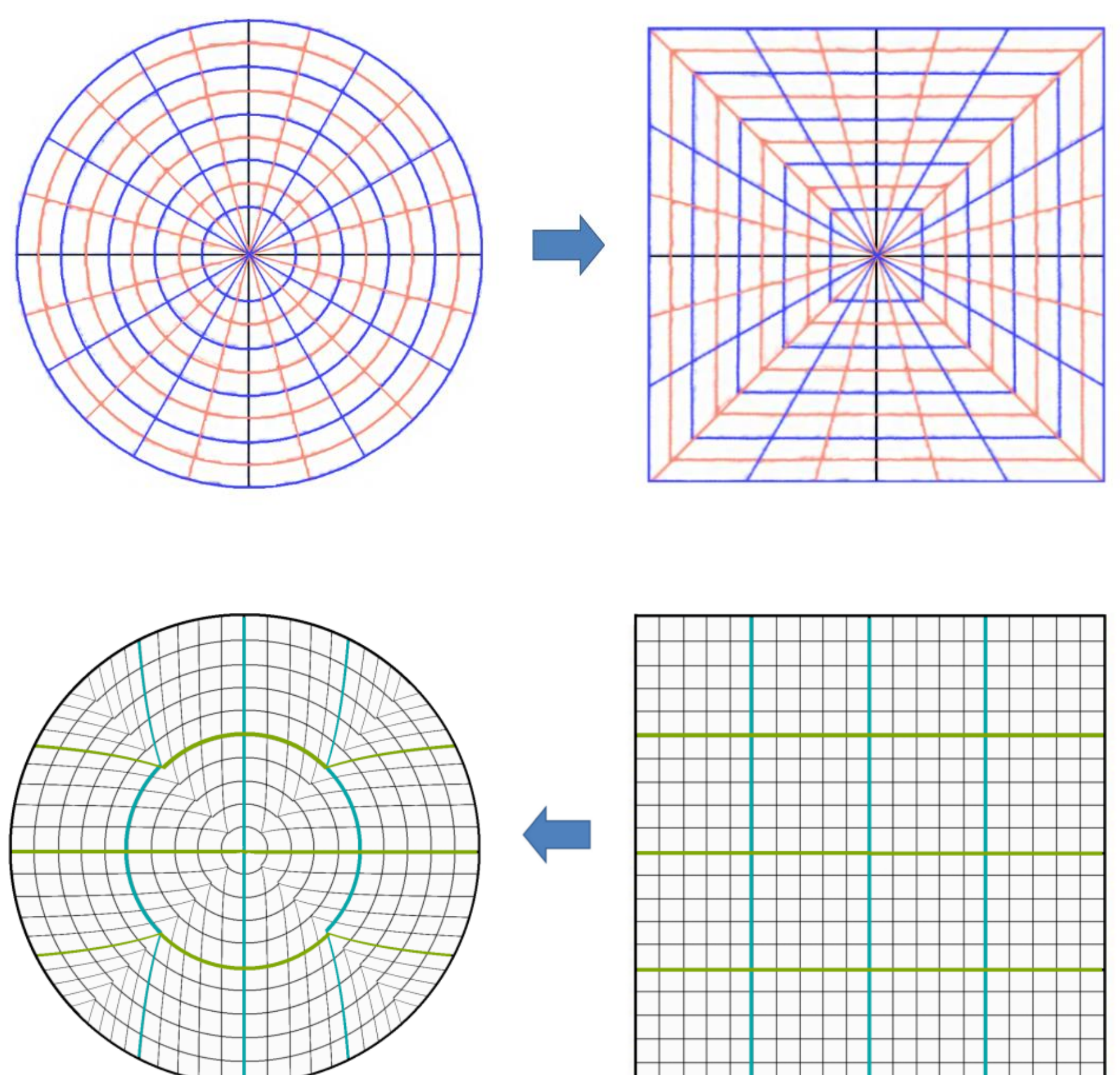

Disc to square mapping:

$$x = \begin{cases} sgn(u)\sqrt{u^2+v^2} & when\ u^2 \ge v^2 \\ sgn(v)\ \dfrac{u}{v}\sqrt{u^2+v^2} & when\ u^2 < v^2 \end{cases}$$

$$y = \begin{cases} sgn(u)\dfrac{v}{u}\ \sqrt{u^2+v^2} & when\ u^2 \ge v^2 \\ sgn(v)\ \sqrt{u^2+v^2} & when\ u^2 < v^2 \end{cases}$$

Square to disc mapping:

$$u = \begin{cases} sgn(x)\dfrac{x^2}{\sqrt{x^2+y^2}} & when\ x^2 \ge y^2 \\ sgn(y)\dfrac{x\,y}{\sqrt{x^2+y^2}} & when\ x^2 < y^2 \end{cases}$$

$$v = \begin{cases} sgn(x)\dfrac{x\,y}{\sqrt{x^2+y^2}} & when\ x^2 \ge y^2 \\ sgn(y)\dfrac{y^2}{\sqrt{x^2+y^2}} & when\ x^2 < y^2 \end{cases}$$

# FG-Squircular Mapping

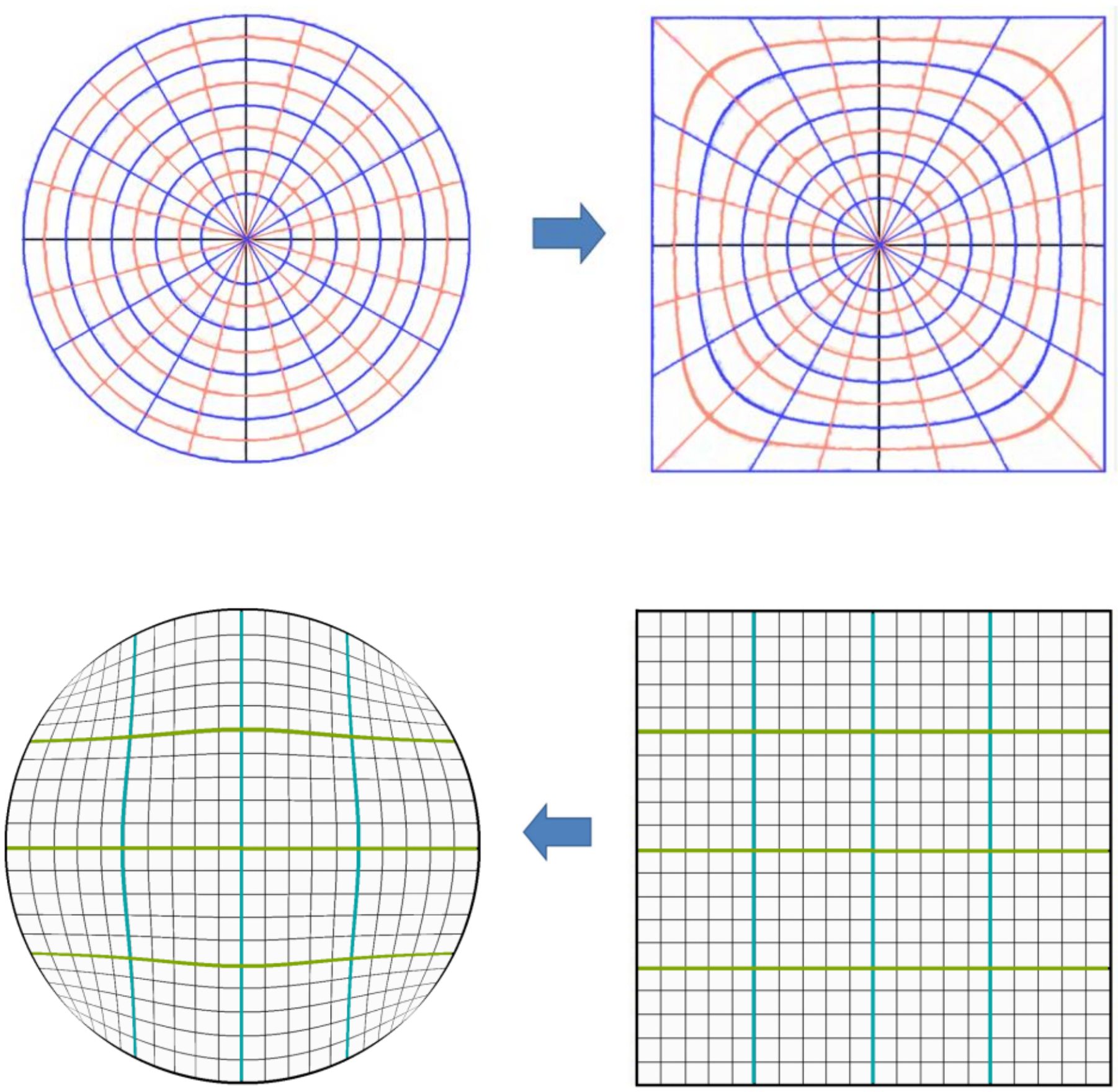

Disc to square mapping:

$$x = \frac{sgn(uv)}{v\sqrt{2}} \sqrt{u^2 + v^2 - \sqrt{(u^2 + v^2)(u^2 + v^2 - 4u^2v^2)}}$$

$$y = \frac{sgn(uv)}{u\sqrt{2}} \sqrt{u^2 + v^2 - \sqrt{(u^2 + v^2)(u^2 + v^2 - 4u^2v^2)}}$$

Square to disc mapping:

$$u = \frac{x\sqrt{x^2 + y^2 - x^2y^2}}{\sqrt{x^2 + y^2}} \qquad v = \frac{y\sqrt{x^2 + y^2 - x^2y^2}}{\sqrt{x^2 + y^2}}$$

# Elliptical Grid Mapping

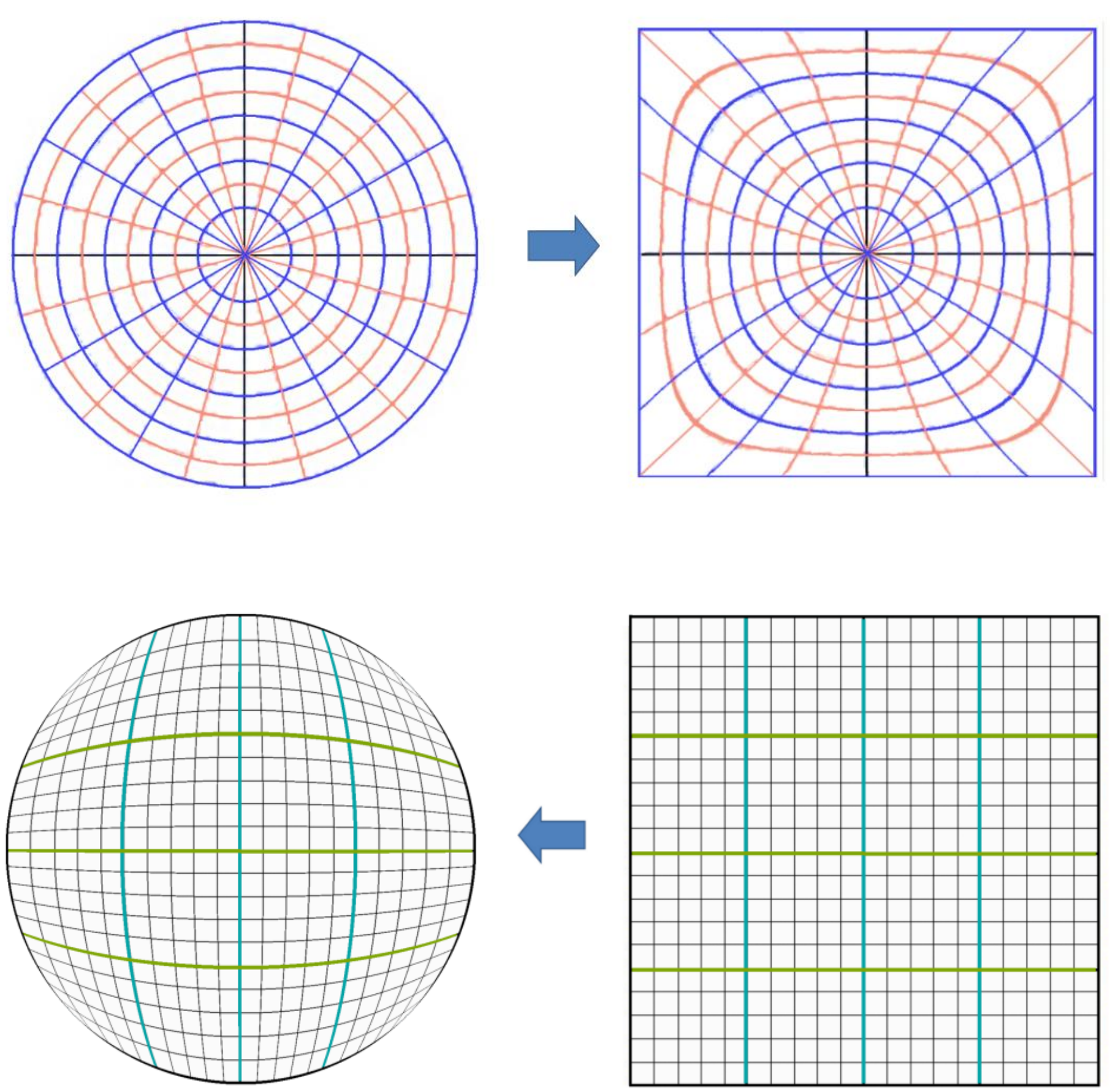

Disc to square mapping:

$$x = \frac{1}{2}\sqrt{2 + u^2 - v^2 + 2\sqrt{2}\,u} - \frac{1}{2}\sqrt{2 + u^2 - v^2 - 2\sqrt{2}\,u}$$

$$y = \frac{1}{2}\sqrt{2 - u^2 + v^2 + 2\sqrt{2}\,v} - \frac{1}{2}\sqrt{2 - u^2 + v^2 - 2\sqrt{2}\,v}$$

Square to disc mapping:

$$u = x\sqrt{1 - \frac{y^2}{2}} \qquad\qquad v = y\sqrt{1 - \frac{x^2}{2}}$$

# Schwarz-Christoffel Mapping

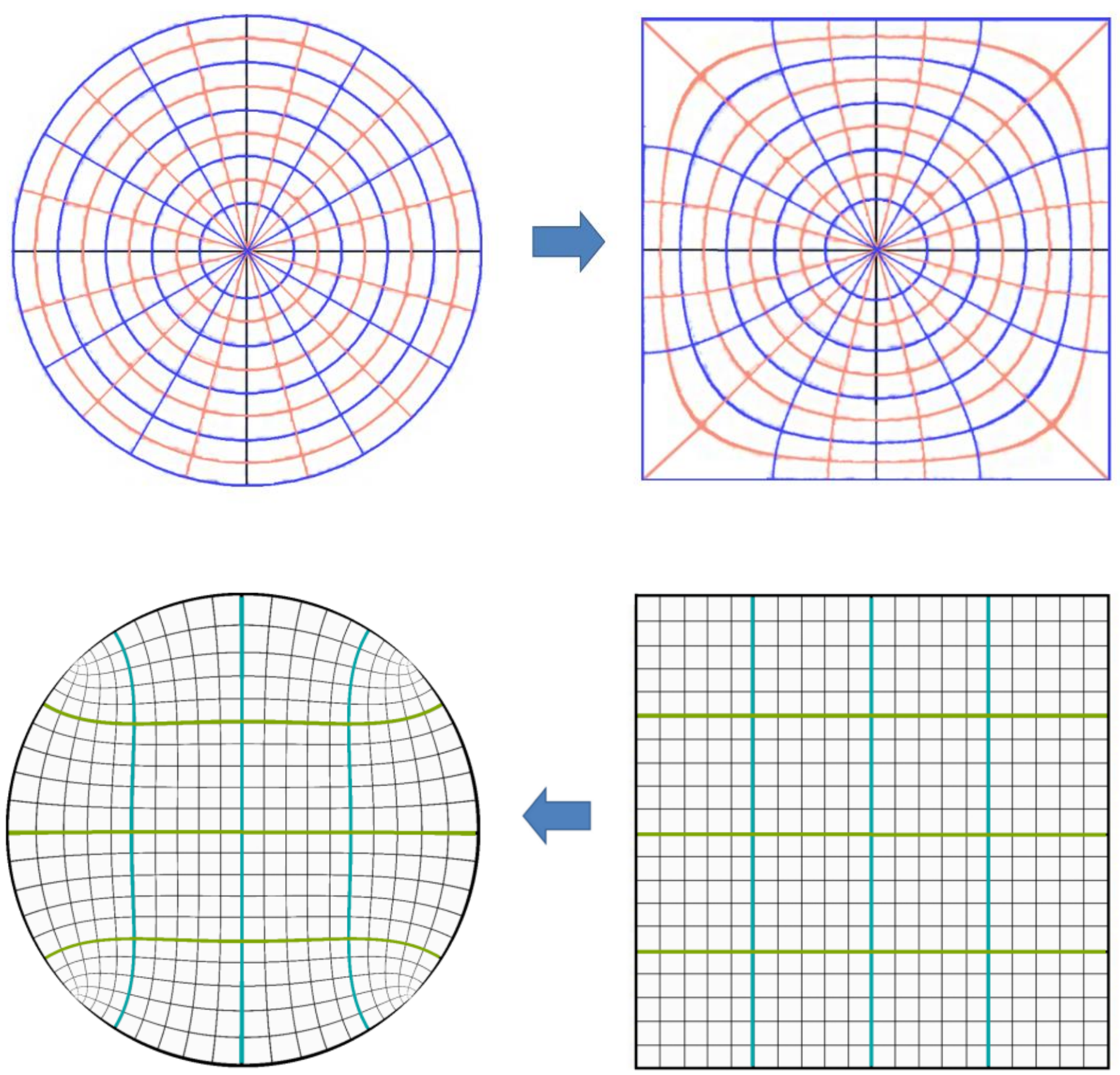

Disc to square mapping:

$$x = Re\left(\frac{1-i}{-K_e} F\left(\cos^{-1}\left(\frac{1+i}{\sqrt{2}}(u+v\,i)\right)\ ,\frac{1}{\sqrt{2}}\right)\right) + 1$$

$$y = Im\left(\frac{1-i}{-K_e} F\left(\cos^{-1}\left(\frac{1+i}{\sqrt{2}}(u+v\,i)\right)\ ,\frac{1}{\sqrt{2}}\right)\right) - 1$$

Square to disc mapping:

$$u = Re\left(\frac{1-i}{\sqrt{2}}\ cn\left(K_e \frac{1+i}{2}(x+y\,i) - K_e\ ,\frac{1}{\sqrt{2}}\right)\right)$$

$$v = Im\left(\frac{1-i}{\sqrt{2}}\ cn\left(K_e \frac{1+i}{2}\ (x+y\,i) - K_e\ ,\frac{1}{\sqrt{2}}\right)\right)$$

where **F** is the incomplete Legendre elliptic integral of the 1st kind
**cn** is a Jacobi elliptic function

$$K_e = \int_0^{\frac{\pi}{2}} \frac{dt}{\sqrt{1-\frac{1}{2}\sin^2 t}} \approx 1.854$$

# Part II. Mathematical Details

## 2. Desirable Properties

### 2.1 Conformal and Equiareal Maps

In our mappings, we want to stay faithful to the source and minimize distortion as much as possible. Two standard metrics used in differential geometry [Feeman 2002][Kuhnel 2006][Floater 2005] are measurements of shape distortion and size distortion. Shape distortion is measured in terms of angle variation between the source and target. Size distortion is measured in terms of area variation between the source and target.

In differential geometry parlance, mappings that preserve angles are called conformal. Similarly, mappings that preserve area are called equiareal [Brown 1935]. The figure below shows an example photograph with contrasting results between a conformal map and equiareal map.

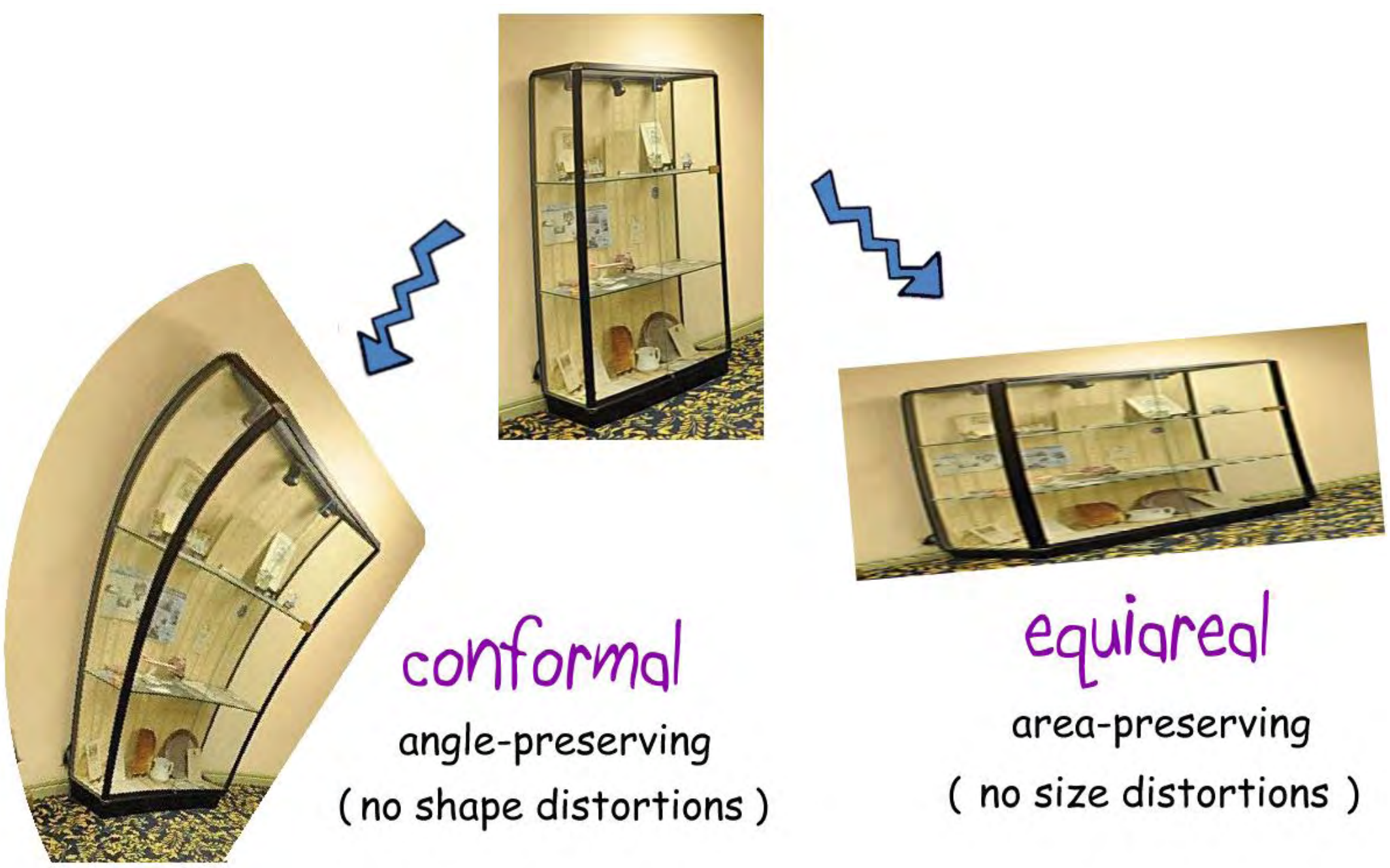

In this paper, we will discuss one conformal mapping -- the Schwarz-Christoffel mapping, which preserves angles throughout the whole mapping. The Schwarz-Christoffel mapping can actually map the unit disc to any simple polygonal region conformally, but we will restrict our discussion to the special case of mapping the circular disc to a square region. The three other mappings mentioned in this paper are not conformal. This can be shown formally by using the Cauchy-Riemann equations in complex analysis.

None of the mappings presented in this paper are equiareal. Although it is possible to alter the mappings given in this paper to make them equiareal, we have not found any nice closed-form analytical expressions for these. Shirley and Chiu did present an equiareal mapping between a circular disc and the square in their 1997 paper [Shirley 1997]. Their mapping has a useful application in ray tracing algorithms for computer graphics [Kolb 1995].

Shirley and Chiu's concentric map has closed-form mapping equations. Qualitatively, the Shirley-Chiu concentric map produces results very similar to the Simple Stretching map, which we will discuss further in this paper. However, unlike the Shirley-Chiu concentric map, the Simple Stretching map is not equiareal.

## 2.2 Radial Mappings

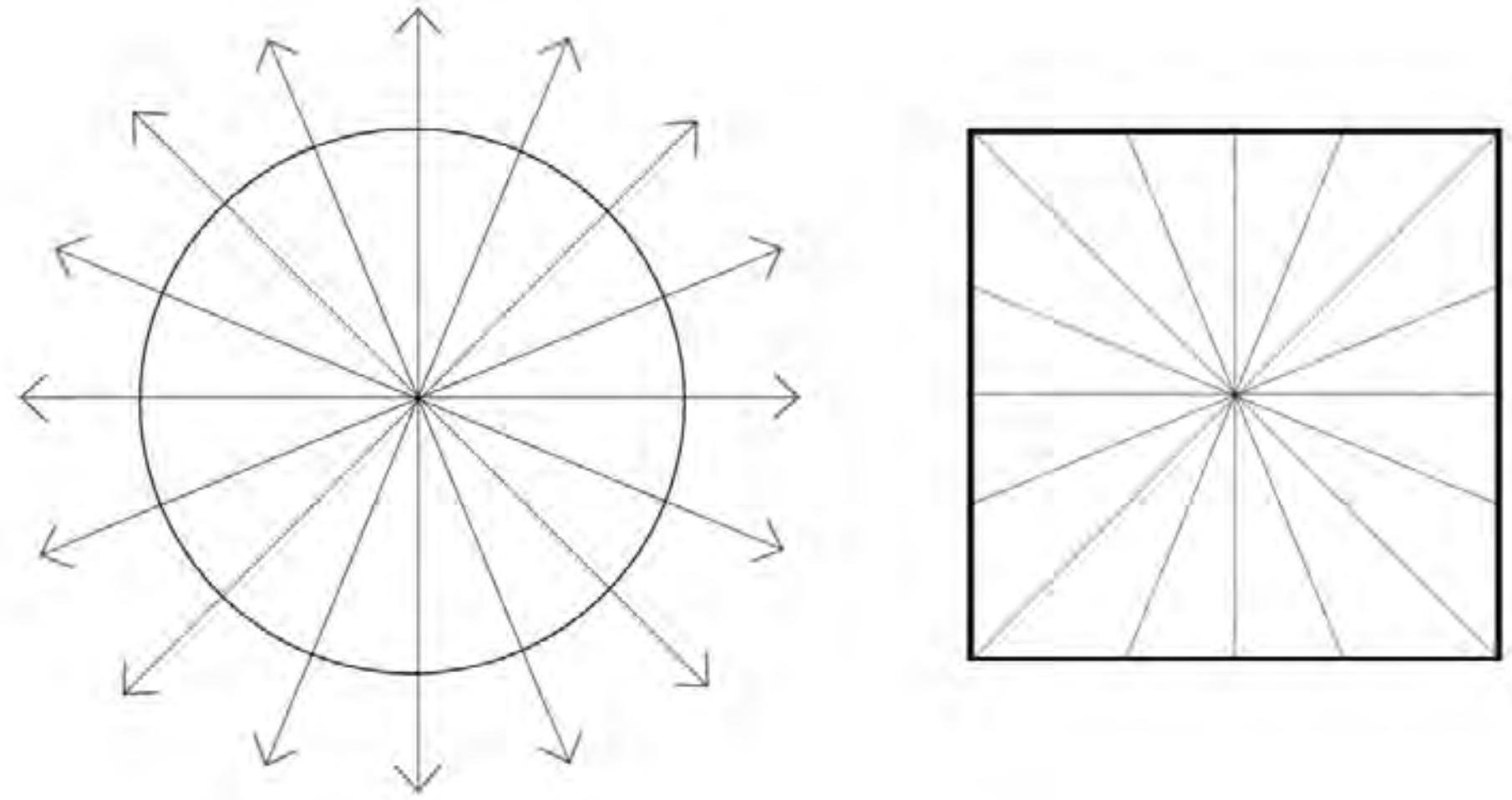

Another desirable property that we want for our mappings is being radial. Intuitively, this means that points can only move along radial lines from the center of the disc during the mapping process. Mathematically, this means that the angle that the point (u,v) makes with the x-axis be the same angle as that of point (x,y).

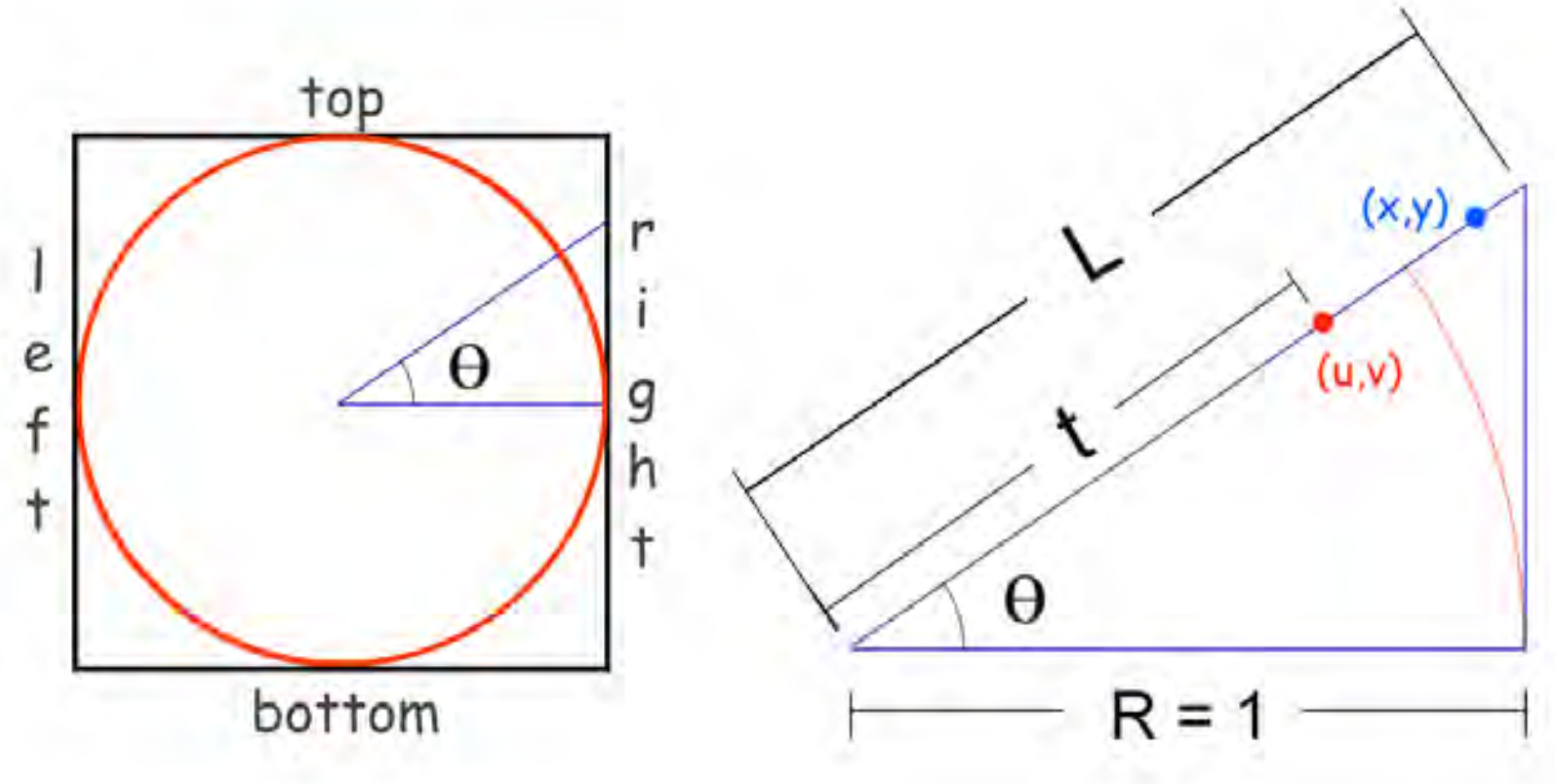

If θ is the angle between the point (u,v) and the x-axis, these equations must hold:

$$\cos\theta = \frac{u}{\sqrt{u^2+v^2}} = \frac{x}{\sqrt{x^2+y^2}}$$

$$\sin\theta = \frac{v}{\sqrt{u^2+v^2}} = \frac{y}{\sqrt{x^2+y^2}}$$

$$\tan\theta = \frac{v}{u} = \frac{y}{x}$$

We call these as the ***radial constraint equations*** for the mapping.

Furthermore, if we write down point (u,v) in polar coordinates (t, θ) inside the unit disc, these equations hold:

$$u = t\ \cos\theta \qquad\qquad v = t \sin\theta$$

Using the radial constraint equations above, we can substitute the trigonometric functions out of the equations to get a relationship between (u,v) and (x,y) for radially constrained mappings.

$$u = t \cos\theta \qquad \Rightarrow \qquad u = t \frac{x}{\sqrt{x^2 + y^2}}$$

$$v = t \sin\theta \qquad \Rightarrow \qquad v = t \frac{y}{\sqrt{x^2 + y^2}}$$

We shall denote these equations as the ***radial mapping linear parametric equations***. Radial mappings between the disc and the square have an equation of this form, with *t* as a properly chosen function of x and y.

## 3. Simple Stretching Map

One of the simplest ways to map a circular disc to a square region is to linearly stretch the circle to the rim of the inscribing square. The equations for stretching from rim to rim are simple but needs to be broken down to four different cases depending on which side of the square the stretching occurs. We consider the case where the circle extends to the right wall. This occurs for angle $\theta$ when $-45° \leq \theta \leq 45°$ . If we parameterize $t$ to be linearly proportional to the distance of the destination point *(x,y)* from the origin, we get:

$$\frac{t}{R} = \frac{\sqrt{x^2 + y^2}}{L}$$

Note that R=1 for our unit circle. Using trigonometry, we have $\cos\theta = \frac{1}{L}$, hence $t = \sqrt{x^2 + y^2} \cos\theta$. Also from trigonometry, we have $\cos\theta = \frac{x}{\sqrt{x^2+y^2}}$ so, the equation simplifies to t = x for the right wall. Using the same reasoning, we can get the value of *t* for the other walls

$$t = \begin{cases} x , & \text{for the right wall} \leftrightarrow \ x \geq |y| \\ y , & \text{for the top wall} \leftrightarrow \ |x| \leq y \\ -x , & \text{for the left wall} \leftrightarrow \ x \leq -|y| \\ -y , & \text{for the bottom wall} \leftrightarrow \ |x| \leq -y \end{cases}$$

Substituting back into the radial mapping parametric equation, we get an equation that relates the point (u,v) in the circular disc to its corresponding point (x,y) in the square.

$$u = \begin{cases} \frac{x^2}{\sqrt{x^2+y^2}}, & \text{for the right wall} \\ \frac{x\,y}{\sqrt{x^2+y^2}}, & \text{for the top wall} \\ \frac{-x^2}{\sqrt{x^2+y^2}}, & \text{for the left wall} \\ \frac{-x\,y}{\sqrt{x^2+y^2}}, & \text{for the bottom wall} \end{cases} \qquad v = \begin{cases} \frac{x\,y}{\sqrt{x^2+y^2}}, & \text{for the right wall} \\ \frac{y^2}{\sqrt{x^2+y^2}}, & \text{for the top wall} \\ \frac{-x\,y}{\sqrt{x^2+y^2}}, & \text{for the left wall} \\ \frac{-y^2}{\sqrt{x^2+y^2}}, & \text{for the bottom wall} \end{cases}$$

Using the signum function and a trick introduced by Dave Cline [Shirley 2011], we can further simplify these equations to:

$$u = \begin{cases} sgn(x) \frac{x^2}{\sqrt{x^2 + y^2}} & when\ x^2 \geq y^2 \\ sgn(y) \frac{x\,y}{\sqrt{x^2 + y^2}} & when\ x^2 < y^2 \end{cases} \qquad v = \begin{cases} sgn(x) \frac{x\,y}{\sqrt{x^2 + y^2}} & when\ x^2 \geq y^2 \\ sgn(y) \frac{y^2}{\sqrt{x^2 + y^2}} & when\ x^2 < y^2 \end{cases}$$

The inverse equations for this Simple Stretching map are simple to derive but require handling 4 different cases for each wall of the square. Start with polar coordinates equation:

$$u = t \ \cos\theta$$

From the radial constraint equations, this becomes

$$u = \frac{t\,u}{\sqrt{u^2+v^2}}$$

which rearranges into this important base equation

$$\boldsymbol{t = \sqrt{u^2+v^2}}$$

**Case 1**: <u>right wall</u> where $t = x$

After substitution, we get $\boldsymbol{x = \sqrt{u^2+v^2}}$ , which is the inverse equation for Case 1 involving *x*.

Now use radial constraint equation: $\frac{v}{u} = \frac{y}{x}$ to get $\frac{v\,y}{u} = \sqrt{u^2+v^2}$

Hence, $\boldsymbol{y = \frac{v}{u}\sqrt{u^2+v^2}}$ is the inverse equation for Case 1 involving *y*.

**Case 2**: <u>top wall</u> where $t = y$

After substitution, we get $\boldsymbol{y = \sqrt{u^2+v^2}}$

Now use radial constraint equation: $\frac{v}{u} = \frac{y}{x}$ to get $\frac{v\,x}{u} = \sqrt{u^2+v^2}$

Hence, $\boldsymbol{x = \frac{u}{v}\sqrt{u^2+v^2}}$ is the inverse equation for Case 2 involving *x*.

**Case 3**: <u>left wall</u> where $t = -x$

After substitution, we get $\boldsymbol{x = -\sqrt{u^2+v^2}}$ , which is the inverse equation for Case 3 involving *x*.

Now use radial constraint equation: $\frac{v}{u} = \frac{y}{x}$ to get $\frac{v\,y}{u} = -\sqrt{u^2+v^2}$

Hence, $\boldsymbol{y = -\frac{v}{u}\sqrt{u^2+v^2}}$ is the inverse equation for Case 3 involving *y*.

**Case 4**: <u>bottom wall</u> where $t = -y$

After substitution, we get $\boldsymbol{y = -\sqrt{u^2+v^2}}$

Now use radial constraint equation: $\frac{v}{u} = \frac{y}{x}$ to get $\frac{v\,x}{u} = -\sqrt{u^2+v^2}$

Hence, $\boldsymbol{x = -\frac{u}{v}\sqrt{u^2+v^2}}$ is the inverse equation for Case 4 involving *x*.

Here is a summary of the inverse equations for each case

$$x = \begin{cases} \sqrt{u^2+v^2}\ , & \text{for the right wall} \\ \frac{u}{v}\sqrt{u^2+v^2}, & \text{for the top wall} \\ -\sqrt{u^2+v^2}\ , & \text{for the left wall} \\ -\frac{u}{v}\sqrt{u^2+v^2}\ , & \text{for the bottom wall} \end{cases}$$

$$y = \begin{cases} \frac{v}{u}\sqrt{u^2+v^2}, & \text{for the right wall} \\ \sqrt{u^2+v^2}\ , & \text{for the top wall} \\ -\frac{v}{u}\sqrt{u^2+v^2}\ , & \text{for the left wall} \\ -\sqrt{u^2+v^2}, & \text{for the bottom wall} \end{cases}$$

We can then use the signum function to simplify these inverse equations. Moreover, we can make further simplifications to the equations if we define and use a ***max*** function:

$$\max(a,b) = \begin{cases} a & when\ a \geq b \\ b & when\ a < b \end{cases}$$

This max function is not a standard function in Mathematics, but it is useful enough that many computer languages incorporate it as a built-in function in their standard libraries. For example, C++ has the max() function in its Standard Template Library. Similarly, Java and Python have a built-in max() function.

In https://marc-b-reynolds.github.io/math/2017/01/08/SquareDisc.html, Marc B. Reynolds used the max function to further simplify the Simple Stretching equations as

Disc to square mapping:

$$\begin{bmatrix} x \\ y \end{bmatrix} = \frac{\sqrt{u^2+v^2}}{\max(|u|,|v|)} \begin{bmatrix} u \\ v \end{bmatrix}$$

Square to disc mapping:

$$\begin{bmatrix} u \\ v \end{bmatrix} = \frac{\max(|x|,|y|)}{\sqrt{x^2+y^2}} \begin{bmatrix} x \\ y \end{bmatrix}$$

Note that the mapping is expressed in terms of a terse vector equation. The variables *x* and *y* are packaged as a 2D vector; likewise for the variables *u* and *v*. Also, the use of the max function effectively nullifies the need for the signum function. This is a welcome change to the equations because the concise form makes it suitable for implementation in various programming languages in the computer.

# 4. Fernandez-Guasti's Squircle

## 4.1 FG-Squircle

In 1992, Manuel Fernandez-Guasti introduced an algebraic equation for representing an intermediate shape between the circle and the square [Fernandez-Guasti 1992]. His equation included a parameter $s$ that can be used to blend the circle and the square smoothly. In this paper, we shall denote this shape as the *Fernandez-Guasti squircle* or *FG-squircle* for short. The figure below illustrates the FG-squircle at varying values of $s$.

The equation for this shape is: $x^2 + y^2 - \frac{s^2}{\rho^2}\, x^2 y^2 = \rho^2$

The squareness parameter $s$ can have any value between 0 and 1. When s = 0, the equation produces a circle with radius $\rho$. When s = 1, the equation produces a square with a side length of $2\rho$. In between, the equation produces a smooth curve that interpolates between the two shapes. Unlike the square, the FG-squircle has no tangent discontinuity along its four corners except at s = 1.

In this paper, we shall restrict our discussion and scope of the FG-squircle to $-\rho \leq x \leq \rho$ and $-\rho \leq y \leq \rho$. The FG-squircle equation is valid in the regions $|x| > \rho$ and $|y| > \rho$, but we will ignore those regions.

## 4.2 The Shrunken FG-Squircle

Consider what happens to the FG-squircle equation when we equate the squareness $s$ with $\rho$. The equation reduces to $x^2 + y^2 - x^2 y^2 = s^2$

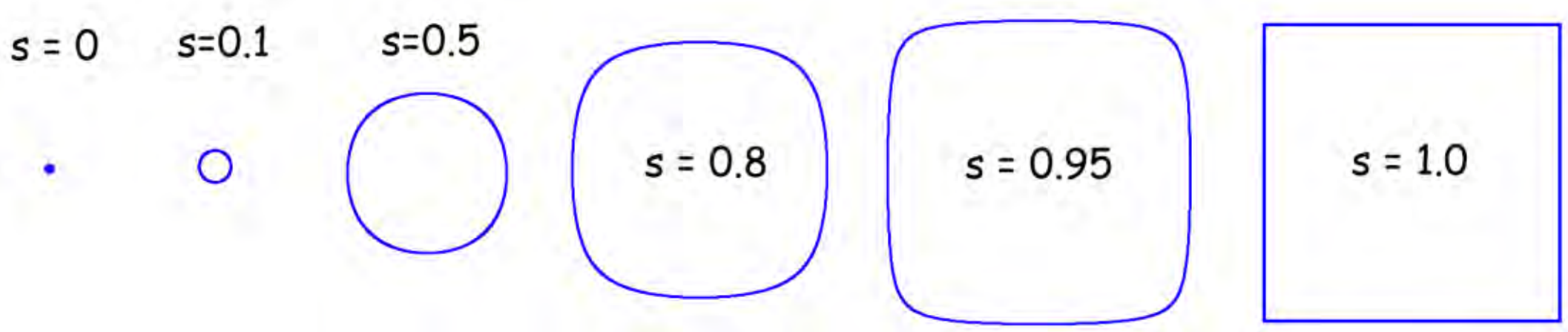

In this equation, we have one parameter $s$ that controls both the squareness and the size of the FG-squircle. Visually, we have a continuum of growing shapes that start as a point at $s=0$ and end as a square at $s=1$. At $s = ½$, we have a half-sized FG-squircle that is not quite circular nor square in shape, but in between the two.

This special case of the FG-squircle when $s = \rho$ is essential to our derivation of several disc-to-square mappings. We shall call this shape as the shrunken FG-squircle.

### 4.3 The Circular Continuum of the Disc

Recall that we defined our input unit disc as the set $\mathcal{D} = \{(u, v) \in \mathbb{R}^2 |\ u^2 + v^2 \leq 1\}$. If we think of the unit disc as a continuum of concentric circles with radii growing from zero to one, we can parameterize the unit disc as the set $\mathcal{D} = \{(u, v) \in \mathbb{R}^2 |\ u^2 + v^2 = t^2 \ \ ,\ \ 0 \leq t \leq 1\}$. In doing so, we introduced a parameter t that is the distance of point (u,v) to the origin.

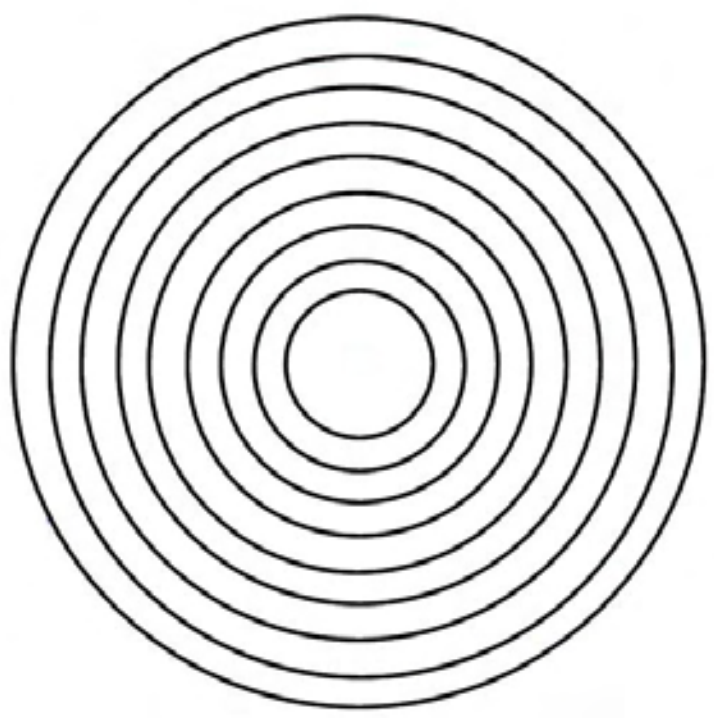

### 4.4 The Squircular Continuum of the Square

In analogy to the circular continuum of the unit disc, one can write the square region [-1,1] x [-1,1] as the set $\mathcal{S} = \{(x, y) \in \mathbb{R}^2 |\ x^2 + y^2 - x^2y^2 = t^2\ ,\ \ 0 \leq t \leq 1\}$

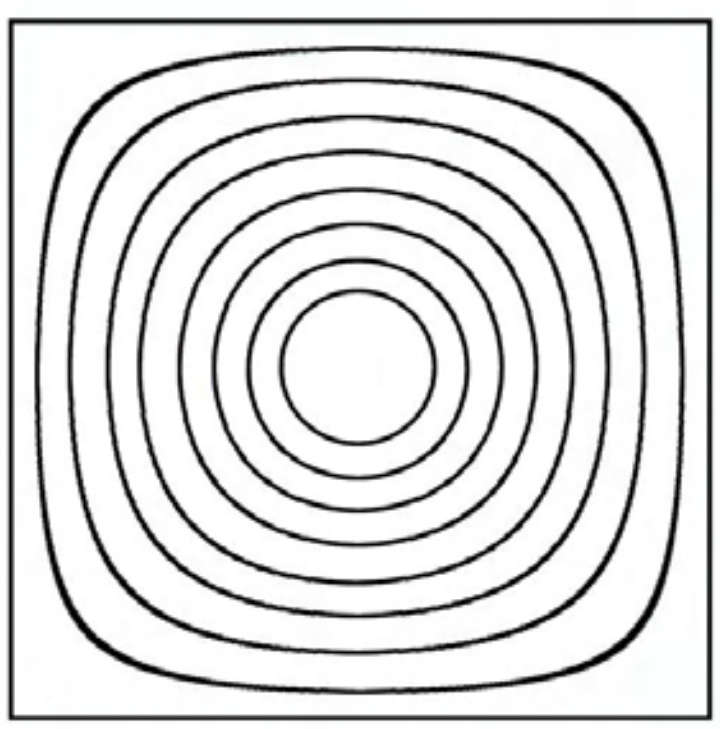

In other words, the square can be considered as a continuum of concentric shrunken FG-squircles.

### 4.5 Linear Squircular Continuum Mapping Equation

Recall from our circular continuum discussion that the unit disc can be represented as a set of concentric circles and parameterized as $\mathcal{D} = \{(u, v) \in \mathbb{R}^2 |\ u^2 + v^2 = t^2 \ \ ,\ \ 0 \leq t \leq 1\}$. Likewise, recall from the previous section that our square region can be represented as a set of concentric shrunken FG-squircles and parameterized as $\mathcal{S} = \{(x, y) \in \mathbb{R}^2 |\ x^2 + y^2 - x^2y^2 = t^2\ ,\ \ 0 \leq t \leq 1\}$.

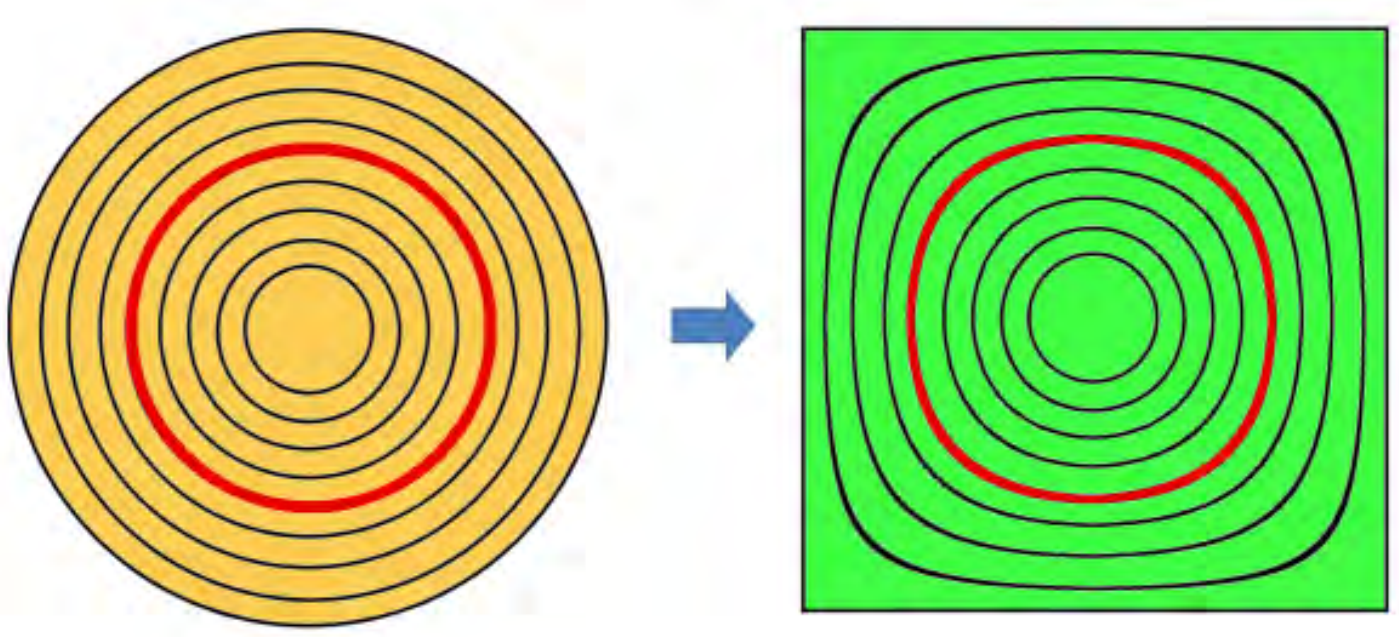

We can establish a correspondence between the unit disc and the square region by mapping every circular contour in the interior of the disc to a squircular contour in the interior of the square. In other words, we map contour curves in the circular continuum of the disc to those in the squircular continuum of the square. This can be done by equating the parameter t of both sets to get the equation:

$$u^2 + v^2 = x^2 + y^2 - x^2 y^2$$

We denote this equation as the *linear squircularity condition* for mapping a circular disc to a square region.

## 5. FG-Squircular Mapping

Using the Fernandez Guasti squircle, we can design a way to map a circular disc smoothly to a square region. The main idea is to map each circular contour in the interior of the disc to a squircular contour in the interior of the square. We can combine this with the radial constraint and derive the mapping equations from there.

### 5.1 Derivation of the Disc to Square Equations

It is easy to derive the FG-Squircular mapping by combining the *linear squircularity condition* from the previous section given by the equation

$$t^2 = u^2 + v^2 = x^2 + y^2 - x^2 y^2$$

with the *radial mapping linear parametric equations*

$$u = t \, \frac{x}{\sqrt{x^2 + y^2}} \qquad\qquad v = t \, \frac{y}{\sqrt{x^2 + y^2}}$$

After substitution of parameter t, we get:

$$u = \frac{x\sqrt{x^2 + y^2 - x^2 y^2}}{\sqrt{x^2 + y^2}} \qquad\qquad v = \frac{y\sqrt{x^2 + y^2 - x^2 y^2}}{\sqrt{x^2 + y^2}}$$

This gives us the FG-Squircular mapping as shown in page 4. This mapping has the nice property of being radial as well as being compliant with the linear squircularity condition. In other words, this is a radial mapping that converts circular contours on the disc to squircular contours on the square.

Only two of the four mappings presented in this paper are radial -- the other being the Simple Stretching map. However, for imaging applications, the FG-Squircular mapping produces significantly better results than the Simple Stretching map. This is because the latter converts circular contours inside the circular disc to square contours inside the square. The square has a $C^1$ discontinuity on its four corners, so the mapping produces $C^1$ discontinuities along the main diagonals

In contrast, the FG-Squircular mapping converts circular contours inside the circular disc to squircular contours inside the square. The FG-squircle does not have tangent discontinuity on its four corners, so there are no diagonal discontinuities in the FG-Squircular mapping.

## 5.2 Inversion of the FG-Squircular Mapping

We shall now derive the inverse equations for the FG-Squircular mapping. Since the FG-Squircular mapping is radial by design, we can use the *radial constraint equations* to get a relationship between x with y.

$$\tan\theta = \frac{v}{u} = \frac{y}{x} \qquad \Rightarrow \qquad y = \frac{v\,x}{u}$$

Substitute this into the equation for the *linear squircularity condition*:

$$u^2 + v^2 = x^2 + y^2 - x^2y^2$$

$$u^2 + v^2 = x^2 + \left(\frac{v\,x}{u}\right)^2 - x^2\left(\frac{v\,x}{u}\right)^2$$

$$= x^2 + \frac{v^2}{u^2}\,x^2 - \frac{v^2}{u^2}x^4$$

$$= \left(1 + \frac{v^2}{u^2}\right)x^2 - x^4\,\frac{v^2}{u^2}$$

Rearranging all the terms to one side of the equation, we get

$$\frac{v^2}{u^2}x^4 - \left(1 + \frac{v^2}{u^2}\right)x^2 + u^2 + v^2 = 0$$

$$\Rightarrow \qquad v^2x^4 - (u^2 + v^2)x^2 + u^4 + u^2v^2 = 0$$

This is a special kind of quartic polynomial equation called a biquadratic. Notice that there are no cubic or linear terms in this 4th degree polynomial equation in x. We can solve for $x^2$ using the quadratic equation with coefficients

$$a = v^2 \qquad\qquad b = -(u^2 + v^2) \qquad\qquad c = u^4 + u^2v^2$$

This gives us the solution for $x^2$

$$x^2 = \frac{u^2 + v^2 \pm \sqrt{(u^2 + v^2)^2 - 4v^2(u^4 + u^2v^2)}}{2v^2}$$

$$= \frac{u^2 + v^2 \pm \sqrt{(u^2 + v^2)^2 - 4u^4v^2 - 4u^2v^4}}{2v^2}$$

$$= \frac{u^2 + v^2 \pm \sqrt{(u^2 + v^2)^2 - 4u^2v^2(u^2 + v^2)}}{2v^2}$$

$$= \frac{u^2 + v^2 \pm \sqrt{(u^2 + v^2)(u^2 + v^2 - 4u^2v^2)}}{2v^2}$$

We can then get a quadrant-aware inverse equation for x as

$$x = \frac{sgn(uv)}{v\sqrt{2}}\sqrt{\;u^2 + v^2 - \sqrt{(u^2 + v^2)(u^2 + v^2 - 4u^2v^2)}}$$

Using the *radial constraint equations* and substituting x, we can also get an inverse equation for y

$$y = \frac{vx}{u} \qquad \Rightarrow \qquad y = \frac{sgn(uv)}{u\sqrt{2}}\sqrt{\;u^2 + v^2 - \sqrt{(u^2 + v^2)(u^2 + v^2 - 4u^2v^2)}}$$

Note that it possible to write these inverse equations in a slightly different form that is more amenable to software implementation in the computer.

We need to use two properties of the signum function. The 1st is a product rule:

$$sgn(uv) = sgn(u)\, sgn(v)$$

Also, we observe that for $v \neq 0$ ,

$$\frac{|v|}{v} = \frac{v}{|v|}$$

so for non-zero inputs, the signum function is its own reciprocal

$$sgn(v) = \frac{1}{sgn(v)} = \frac{v}{|v|}$$

$$\therefore \qquad \frac{sgn(v)}{v} = \frac{1}{|v|}$$

This 2nd property follows easily from the definition of the signum function.

So, we can rewrite the equations for x as:

$$x = \frac{sgn(uv)}{v\sqrt{2}} \sqrt{u^2 + v^2 - \sqrt{(u^2 + v^2)(u^2 + v^2 - 4u^2v^2)}}$$

$$= \frac{sgn(u)\, sgn(v)}{v\sqrt{2}} \sqrt{u^2 + v^2 - \sqrt{(u^2 + v^2)(u^2 + v^2 - 4u^2v^2)}}$$

$$= \frac{sgn(u)}{|v|\sqrt{2}} \sqrt{u^2 + v^2 - \sqrt{(u^2 + v^2)(u^2 + v^2 - 4u^2v^2)}}$$

Likewise, we can get the equation for y as:

$$y = \frac{sgn(v)}{|u|\sqrt{2}} \sqrt{u^2 + v^2 - \sqrt{(u^2 + v^2)(u^2 + v^2 - 4u^2v^2)}}$$

# 6. Elliptical Grid Mapping

In 2005, Philip Nowell introduced a square to disc mapping that converts horizontal and vertical lines in the square to elliptical arcs inside a circular region. This mapping turns a regular rectangular grid into a regular curvilinear grid consisting of elliptical arcs. Nowell provided a mathematical derivation of his mapping in his blog [Nowell 2005] which we will not repeat here. However, he left out the reversal of the process. In this section, we shall provide inverse equations to his mapping.

## 6.1 Inversion using Trigonometry

We now derive the inverse of Nowell's Elliptical Grid mapping. The derivation involves 3 steps which is summarized as

Step 1: Convert circular coordinates (u,v) to polar coordinates and introduce trigonometric variables α and β

Step 2: Find an expression for (x,y) in terms of α and β

Step 3: Find an expression for (x,y) in terms of u and v

Start with Nowell's equations:

$$u = x\sqrt{1 - \frac{y^2}{2}} \qquad\qquad v = y\sqrt{1 - \frac{x^2}{2}}$$

*Step 1: Convert circle coordinates* (u,v) *to polar form* (r, θ)

$$r = \sqrt{u^2 + v^2}$$

$$tan\,\theta = \frac{v}{u} \qquad sin\,\theta = \frac{v}{\sqrt{u^2 + v^2}} \qquad cos\,\theta = \frac{u}{\sqrt{u^2 + v^2}}$$

We introduce intermediate trigonometric angles α and β such that

$$\cos\alpha = r\ \cos(\theta + \frac{\pi}{4})$$
$$\cos\beta = r\ \cos(\theta - \frac{\pi}{4})$$

Note that, since r≤ 1 and the cosine value is ≤ 1, there exist angles α and β that satisfies this. We can then expand these expressions using trigonometric expansion formulas.

$$\cos\alpha = r\ \cos(\theta + \frac{\pi}{4}) = r\left(\cos\theta\cos\frac{\pi}{4} - \sin\theta\sin\frac{\pi}{4}\right) = r\frac{\sqrt{2}}{2}(\cos\theta - \sin\theta)$$

Likewise,

$$\cos\beta = r\ \cos(\theta - \frac{\pi}{4}) = r\left(\cos\theta\cos\frac{\pi}{4} + \sin\theta\sin\frac{\pi}{4}\right) = r\frac{\sqrt{2}}{2}(\cos\theta + \sin\theta)$$

Furthermore, we can find expressions for cos α and cos β in terms of u and v.

$$\cos\alpha = r\frac{\sqrt{2}}{2}(\cos\theta - \sin\theta) = r\frac{\sqrt{2}}{2}\left(\frac{u}{\sqrt{u^2+v^2}} - \frac{v}{\sqrt{u^2+v^2}}\right) = \sqrt{u^2+v^2}\,\frac{\sqrt{2}}{2}\frac{u-v}{\sqrt{u^2+v^2}} = \frac{\sqrt{2}}{2}(u-v)$$

$$\cos\beta = r\frac{\sqrt{2}}{2}(\cos\theta + \sin\theta) = r\frac{\sqrt{2}}{2}\left(\frac{u}{\sqrt{u^2+v^2}} + \frac{v}{\sqrt{u^2+v^2}}\right) = \sqrt{u^2+v^2}\,\frac{\sqrt{2}}{2}\frac{u+v}{\sqrt{u^2+v^2}} = \frac{\sqrt{2}}{2}(u+v)$$

$$\therefore \qquad \cos\alpha = \frac{\sqrt{2}}{2}(u-v) \qquad\qquad \cos\beta = \frac{\sqrt{2}}{2}(u+v)$$

*Step 2: Find an expression for x and y in terms of angle $\alpha$ and $\beta$*

From step 1, we can write cos α and cos β in terms of x and y

$$\cos\alpha = \frac{\sqrt{2}}{2}\left(\mathrm{x}\sqrt{1-\frac{y^2}{2}} - \mathrm{y}\sqrt{1-\frac{x^2}{2}}\right) = \frac{x}{\sqrt{2}}\sqrt{1-\frac{y^2}{2}} - \frac{y}{\sqrt{2}}\sqrt{1-\frac{x^2}{2}}$$

$$\cos\beta = \frac{\sqrt{2}}{2}\left(\mathrm{x}\sqrt{1-\frac{y^2}{2}} + \mathrm{y}\sqrt{1-\frac{x^2}{2}}\right) = \frac{x}{\sqrt{2}}\sqrt{1-\frac{y^2}{2}} + \frac{y}{\sqrt{2}}\sqrt{1-\frac{x^2}{2}}$$

We now introduce another set of intermediate trigonometric variables ϕ and λ. Define:

$$\phi = \cos^{-1}\frac{x}{\sqrt{2}}$$

$$\lambda = \sin^{-1}\frac{y}{\sqrt{2}}$$

Note that since x≤1 and y≤1, it is okay to compute their inverse trigonometric values. This implies

$$\cos\phi = \frac{x}{\sqrt{2}} \qquad \sin\phi = \sqrt{1-\frac{x^2}{2}}$$

and

$$\cos\lambda = \sqrt{1-\frac{y^2}{2}} \qquad \sin\lambda = \frac{y}{\sqrt{2}}$$

So we can write cos α and cos β in terms of ϕ and λ

$$\cos\alpha = \cos\phi\cos\lambda - \sin\phi\sin\lambda = \cos(\phi+\lambda)$$

$$\cos\beta = \cos\phi\cos\lambda + \sin\phi\sin\lambda = \cos(\phi-\lambda)$$

Hence

$$\alpha = \phi + \lambda \ = \cos^{-1}\frac{x}{\sqrt{2}} + \sin^{-1}\frac{y}{\sqrt{2}}$$

$$\beta = \phi - \lambda = \cos^{-1}\frac{x}{\sqrt{2}} - \sin^{-1}\frac{y}{\sqrt{2}}$$

Taking the sum and difference of these two equations, we get

$$\alpha + \beta = 2\cos^{-1}\frac{x}{\sqrt{2}}$$

$$\alpha - \beta = 2\sin^{-1}\frac{y}{\sqrt{2}}$$

Rearranging the terms, we get

$$x = \sqrt{2}\cos(\frac{\alpha+\beta}{2})$$

$$y = \sqrt{2}\sin(\frac{\alpha-\beta}{2})$$

*Step 3: get an expression for x and y in terms of u and v*

We start with the result from step 2 and expand using trigonometric identities for sums and half-angles

$$x = \sqrt{2}\ \cos(\frac{\alpha}{2}+\frac{\beta}{2}) = \ \sqrt{2}(\cos\frac{\alpha}{2}\cos\frac{\beta}{2} - \sin\frac{\alpha}{2}\sin\frac{\beta}{2}) = \ \sqrt{2}(\sqrt{\frac{(1+\cos\alpha)}{2}}\sqrt{\frac{(1+\cos\beta)}{2}} - \sqrt{\frac{(1-\cos\alpha)}{2}}\sqrt{\frac{(1-\cos\beta)}{2}}\ )$$

$$= \frac{1}{2}\sqrt{2(1+\cos\alpha)(1+\cos\beta)} \ - \frac{1}{2}\sqrt{2(1-\cos\alpha)(1-\cos\beta)}$$

$$y = \ \sqrt{2}\sin(\frac{\alpha}{2}-\frac{\beta}{2}) = \sqrt{2}(\sin\frac{\alpha}{2}\cos\frac{\beta}{2} - \cos\frac{\alpha}{2}\sin\frac{\beta}{2}) = \ \sqrt{2}\left(\sqrt{\frac{(1-\cos\alpha)}{2}}\sqrt{\frac{(1+\cos\beta)}{2}} - \sqrt{\frac{(1+\cos\alpha)}{2}}\sqrt{\frac{(1-\cos\beta)}{2}}\right)$$

$$= \frac{1}{2}\sqrt{2(1+\cos\alpha)(1+\cos\beta)} \ - \frac{1}{2}\sqrt{2(1-\cos\alpha)(1-\cos\beta)}$$

Now, recall from step 1 that

$$\cos\alpha = \frac{\sqrt{2}}{2}(u - v) \qquad\qquad \cos\beta = \frac{\sqrt{2}}{2}(u + v)$$

So we can substitute u and v values into cos α and cos β

$$x = \frac{1}{2}\sqrt{2(1+\cos\alpha)(1+\cos\beta)} - \frac{1}{2}\sqrt{2(1-\cos\alpha)(1-\cos\beta)}$$

$$= \frac{1}{2}\sqrt{2\left(1+\frac{\sqrt{2}}{2}(u-v)\right)\left(1+\frac{\sqrt{2}}{2}(u+v)\right)} - \frac{1}{2}\sqrt{2\left(1-\frac{\sqrt{2}}{2}(u-v)\right)\left(1-\frac{\sqrt{2}}{2}(u+v)\right)}$$

$$\boldsymbol{= \frac{1}{2}\sqrt{2+u^2-v^2+2\sqrt{2}\,u} - \frac{1}{2}\sqrt{2+u^2-v^2-2\sqrt{2}\,u}}$$

Likewise,

$$y = \frac{1}{2}\sqrt{2(1-\cos\alpha)(1+\cos\beta)} - \frac{1}{2}\sqrt{2(1+\cos\alpha)(1-\cos\beta)}$$

$$= \frac{1}{2}\sqrt{2\left(1-\frac{\sqrt{2}}{2}(u-v)\right)\left(1+\frac{\sqrt{2}}{2}(u+v)\right)} - \frac{1}{2}\sqrt{2\left(1+\frac{\sqrt{2}}{2}(u-v)\right)\left(1-\frac{\sqrt{2}}{2}(u+v)\right)}$$

$$\boldsymbol{= \frac{1}{2}\sqrt{2-u^2+v^2+2\sqrt{2}\,v} - \frac{1}{2}\sqrt{2-u^2+v^2-2\sqrt{2}v}}$$

This completes the derivation of denested inverse equations.

## 6.2 Inversion using the Biquadratic Equation

The inverse equations for the Elliptical Grid mapping we derived in the previous section are by no means unique in form. We can actually derive another set of inverse equations by using a different method. Of course, these sets of inverse equations are ultimately equivalent to each other. That is, they are just different manifestations of the same inverse equations.

The inverse equations that we will derive in this section are not as mathematically elegant as the denested equations previously derived. Nevertheless, the inverse equations here are valid and just as useful.

We look at Nowell's square to disc equations and derive another set of the inverse equations for it. Start with

$$u = \mathrm{x}\sqrt{1-\frac{y^2}{2}} \qquad\qquad v = \mathrm{y}\sqrt{1-\frac{x^2}{2}}$$

Isolate y in the 2nd equation to get:

$$y = \frac{v}{\sqrt{1-\frac{x^2}{2}}} = \frac{\sqrt{2}v}{\sqrt{2-x^2}}$$

$$\Rightarrow \qquad \mathrm{y}^2 = \frac{2\mathrm{v}^2}{2-\mathrm{x}^2}$$

Substituting back to the 1st equation

$$u = x\sqrt{1-\frac{2v^2}{2(2-x^2)}} = x\sqrt{\frac{2-x^2-v^2}{2-x^2}}$$

$$\Rightarrow \qquad u^2 = \frac{x^2(2 - x^2 - v^2)}{2 - x^2}$$

$$\Rightarrow \quad (2 - x^2)u^2 = x^2(2 - x^2 - v^2)$$

$$\Rightarrow \quad x^4 - x^2(2 + u^2 - v^2) + 2u^2 = 0$$

This is a special kind of quartic polynomial equation called a biquadratic. Notice that there are no cubic or linear terms in the 4th degree polynomial equation in x. We can solve for $x^2$ using the quadratic equation with coefficients

$$a = 1 \qquad\qquad b = -(2 + u^2 - v^2) \qquad\qquad c = 2u^2$$

This gives us the solution:

$$x^2 = \frac{2 + u^2 - v^2 \pm \sqrt{(2 + u^2 - v^2)^2 - 8u^2}}{2}$$

We can then get a quadrant-aware inverse equation for x as

$$\boldsymbol{x = \frac{sgn(u)}{\sqrt{2}} \sqrt{2 + u^2 - v^2 - \sqrt{(2 + u^2 - v^2)^2 - 8u^2}}}$$

Using a similar approach, one get the equation for y as

$$\boldsymbol{y = \frac{sgn(v)}{\sqrt{2}} \sqrt{2 - u^2 + v^2 - \sqrt{(2 - u^2 + v^2)^2 - 8v^2}}}$$

### 6.3 Squircularity of the Mapping

It is appropriate at this point to ask what sort of shape does the Elliptical Grid mapping convert circular contours into inside the square. First, recall the circular continuum representation of the unit disc. Now we want to find out what sort of curve inside the square corresponds to each of the concentric circles inside the disc.

Actually, we will show that central circles in the unit disc are mapped into Fernandez-Guasti squircles. We do this by looking at an individual circle inside the disc and substitute x and y values into it.

$$\begin{aligned} u^2 + v^2 &= \left( x\sqrt{1 - \frac{y^2}{2}} \right)^2 + \left( y\sqrt{1 - \frac{x^2}{2}} \right)^2 \\ &= x^2\left(1 - \frac{y^2}{2}\right) + y^2\left(1 - \frac{x^2}{2}\right) \\ &= x^2 - \frac{1}{2}x^2y^2 + y^2 - \frac{1}{2}x^2y^2 \\ &= x^2 + y^2 - x^2y^2 \end{aligned}$$

This is just the equation for the *linear squircularity condition* of the mapping as previously mentioned. From this, we can deduce that the Elliptical Grid mapping also converts central circles in the unit disc to Fernandez-Guasti squircles on the square. *In other words, the Fernandez-Guasti squircle also plays a central role in this mapping!*

### 6.4 Numerical Imprecision in Floating Point Implementations

Recall that the square to disc equations for this mapping are

$$x = \frac{1}{2}\sqrt{2 + u^2 - v^2 + 2\sqrt{2}\,u} - \frac{1}{2}\sqrt{2 + u^2 - v^2 - 2\sqrt{2}\,u}$$
$$y = \frac{1}{2}\sqrt{2 - u^2 + v^2 + 2\sqrt{2}\,v} - \frac{1}{2}\sqrt{2 - u^2 + v^2 - 2\sqrt{2}\,v}$$

In theory, the expressions involving $u$ and $v$ inside the square root operator should never be negative for all *(u, v)* inside the unit disc. In practice though, there is numerical imprecision in floating point arithmetic. This imprecision could result in a very small negative term inside the square roots. This in turn will cause unwanted errors or exceptions within a computer program.

For example, consider the point *(u, v)* = $(\sqrt{2}/2, \sqrt{2}/2)$. The expression $2 + u^2 - v^2 - 2\sqrt{2}\,u$ simplifies to the this expression

$$2 - 2\sqrt{2} * \sqrt{2}/2$$

In some floating point implementations, this expression produces the value of -4.4409e-016, which is a very small negative number. Applying the square root operator on this negative result will cause havoc in some computer programs. This problem occurs in this simple Python program.

```
import math

u = math.sqrt(2)/2
v = math.sqrt(2)/2
result = 2 + u*u – v*v – 2*math.sqrt(2)*u
print(result)
```

Therefore, it is beneficial to add an absolute value function to all the terms inside the square roots in the mapping equations. We can then amend the square to disc equations into

$$\boldsymbol{x = \frac{1}{2}\sqrt{\left|2 + u^2 - v^2 + 2\sqrt{2}\,u\right|} - \frac{1}{2}\sqrt{\left|2 + u^2 - v^2 - 2\sqrt{2}\,u\right|}}$$
$$\boldsymbol{y = \frac{1}{2}\sqrt{\left|2 - u^2 + v^2 + 2\sqrt{2}\,v\right|} - \frac{1}{2}\sqrt{\left|2 - u^2 + v^2 - 2\sqrt{2}\,v\right|}}$$

or equivalently, this mapping equations with nested square roots

$$\boldsymbol{x = \frac{sgn(u)}{\sqrt{2}}\sqrt{2 + u^2 - v^2 - \sqrt{\left|(2 + u^2 - v^2)^2 - 8u^2\right|}}}$$

$$\boldsymbol{y = \frac{sgn(v)}{\sqrt{2}}\sqrt{2 - u^2 + v^2 - \sqrt{\left|(2 - u^2 + v^2)^2 - 8v^2\right|}}}$$

## 7. Schwarz-Christoffel Mapping

In the 1860s, Hermann A. Schwarz and Elwin Christoffel used complex analysis to independently develop a conformal map of the circular disc onto a simple polygonal region. In this paper, we are mainly interested in the special case involving the square. Our goal is to reduce the Schwarz-Christoffel mapping to our specific case and derive an equation relating point *(u,v)* inside the circular disc to point *(x,y)* inside the square.

### 7.1 Fundamental Mapping in the Complex Plane

Without getting into the nitty-gritty details of the Schwarz-Christoffel mapping for now, we show in the diagram below a fundamental conformal mapping between the circular disc and the square in the complex plane. Using the complex-valued Jacobi elliptic function *cn* , one can map every point inside the unit disc to a square region conformally. We wil fully derive this diagram using basic principles of the Schwarz-Christoffel mapping in the next section of this paper. The inverse of the mapping requires the complex-valued incomplete Legendre elliptic integral of the 1st kind *F*.

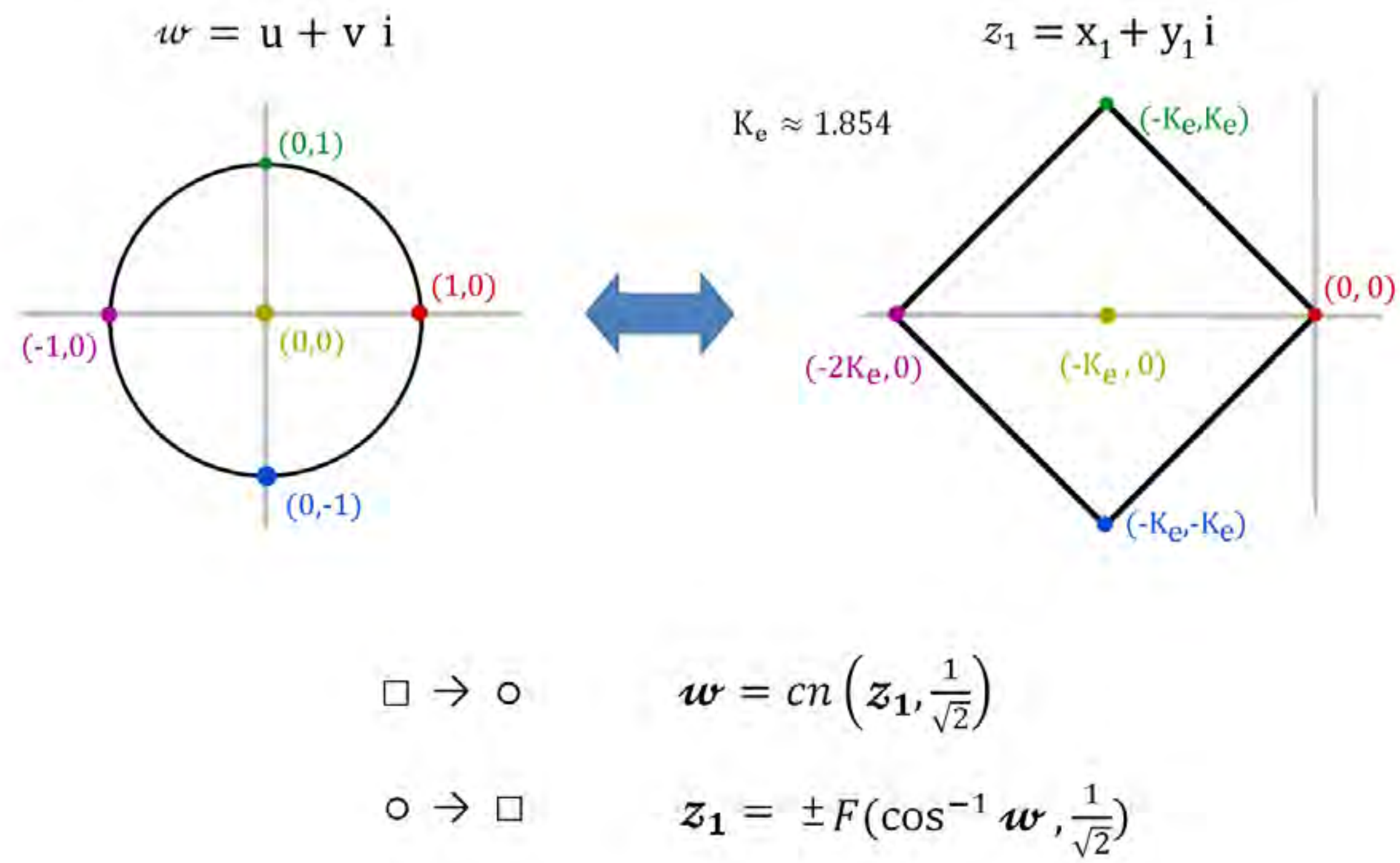

The main drawback of this diagram on the complex plane is that x and y coordinates are not in our canonical mapping space. As the figure above shows, the square has corner coordinates in terms of a constant $K_e$ instead of the ±1 that we desire. Moreover, the square is tilted by 45° and off-center from the origin.

The constant $K_e$ has a numerical value of approximately 1.854. Its exact value is the complete Legendre elliptic integral of the 1st kind with modulus $\frac{1}{\sqrt{2}}$. This number arises from Schwarz and Christoffel's equations for the specific case when the desired polygonal shape is a square.

$$K_e = \int_0^{\frac{\pi}{2}} \frac{dt}{\sqrt{1 - \frac{1}{2}\sin^2 t}} \approx 1.854$$

Note that the complete Legendre elliptic integral of the 1st kind *K* can be calculated from the incomplete Legendre elliptic integral of the 1st kind *F*; i.e.

$$K_e = K\left(\frac{1}{\sqrt{2}}\right) = F\left(\frac{\pi}{2}, \frac{1}{\sqrt{2}}\right)$$

### 7.2 Aligning to the Canonical Mapping Space

In order to get the mapping into our canonical mapping space, we have to do a series of affine transformations on the square. Specifically, we want the square to be centered on the origin with side length value of 2. This is exactly what the square to disc mapping equations below do.

$$u = Re\left(\frac{1-i}{\sqrt{2}}\ cn\left(K_e \frac{1+i}{2}(x+y\,i) - K_e\ ,\frac{1}{\sqrt{2}}\right)\right)$$
$$v = Im\left(\frac{1-i}{\sqrt{2}}\ cn\left(K_e \frac{1+i}{2}\ (x+y\,i) - K_e\ ,\frac{1}{\sqrt{2}}\right)\right)$$

In the complex plane, rotation can be done simply by multiplication with the complex number $e^{i\theta}$. For our case, we are interested in rotation by ±45°, so we have these multipliers

for 45° rotation:

$$e^{i\frac{\pi}{4}} = \cos\frac{\pi}{4} + i\ \sin\frac{\pi}{4} = \frac{1}{\sqrt{2}}(1+i)$$

for -45° rotation:

$$e^{-i\frac{\pi}{4}} = \cos\frac{\pi}{4} - i\ \sin\frac{\pi}{4} = \frac{1}{\sqrt{2}}(1-i)$$

One can locate these rotational multipliers along with $K_e$ translational offsets and scale factors in the mapping equations. These affine transformations are there to place the square into our canonical mapping space in the complex plane. Also, take note that $45° = \pi/4$ radians.

In order to get the disc to square mapping inverse equation, we use the fact that inverse of the Jacobi elliptic function *cn* is just the Legendre elliptic integral of the 1st kind for an arccosine. In other words,

$$cn^{-1}\left(z,\frac{1}{\sqrt{2}}\right) = F(\cos^{-1} z,\frac{1}{\sqrt{2}})$$

Consequently, by using this identity and doing some algebra, the disc to square equations are:

$$x = Re\left(\frac{1-i}{-K_e} F\left(\cos^{-1}\left(\frac{1+i}{\sqrt{2}}(u+v\,i)\right)\ ,\frac{1}{\sqrt{2}}\right)\right) + 1$$

$$y = Im\left(\frac{1-i}{-K_e} F\left(\cos^{-1}\left(\frac{1+i}{\sqrt{2}}(u+v\,i)\right)\ ,\frac{1}{\sqrt{2}}\right)\right) - 1$$

We would like to point out that the complex-valued Jacobi elliptic function *cn* is an even function, so $cn(z,\frac{1}{\sqrt{2}}) = cn\left(-z,\frac{1}{\sqrt{2}}\right)$. Therefore, the equations below are also valid formulas for the mapping:

$$u = Re\left(\frac{1-i}{\sqrt{2}}\ cn\left(-K_e \frac{1+i}{2}(x+y\,i) + K_e\ ,\frac{1}{\sqrt{2}}\right)\right)$$
$$v = Im\left(\frac{1-i}{\sqrt{2}}\ cn\left(-K_e \frac{1+i}{2}\ (x+y\,i) + K_e\ ,\frac{1}{\sqrt{2}}\right)\right)$$

As a consequence of the Jacobi elliptic function *cn* being an even function, its inverse equation can be written more appropriately as

$$cn^{-1}\left(z, \frac{1}{\sqrt{2}}\right) \;=\; \pm F(\cos^{-1} z, \frac{1}{\sqrt{2}})$$

This explains the appearance of the ± sign in the fundamental conformal diagram shown earlier.

## 7.3 A Compact Complex Equation

It is possible to get a more compact complex equation for the mapping on the complex plane. First, observe that the 45° rotational factor needed for canonical alignment can expressed as the square root of *i*. In other words,

$$\frac{1}{\sqrt{2}}(1+i) \;=\; \sqrt{i}$$

Likewise, the -45° rotational factor can be expressed as the square root of *–i*.

$$\frac{1}{\sqrt{2}}(1-i) = \sqrt{-i}$$

Note that these two rotational factors are complex conjugates of each other.

If we define complex variables $\boldsymbol{w}$ and $\boldsymbol{z}$ as

$$\boldsymbol{w} = u + v\,i \qquad\qquad \boldsymbol{z} = x + y\,i$$

i.e. $\boldsymbol{w}$ represents complex numbers inside the circular disc and $\boldsymbol{z}$ represents complex numbers inside the square.

We can then simplify the canonically-aligned mapping equation into

$$\boldsymbol{w} = \sqrt{-i}\;\; cn\left(K_e \boldsymbol{z}\sqrt{\frac{i}{2}} - K_e, \frac{1}{\sqrt{2}}\right)$$

along with its inverse

$$\boldsymbol{z} = \; 1 - i - \frac{\sqrt{-2i}}{K_e}\; F\left(\cos^{-1}\!\left(\boldsymbol{w}\sqrt{i}\right), \frac{1}{\sqrt{2}}\right)$$

These are more compact equations than the ones provided in page 6. Moreover, expressing the equations in this form makes it is easy to see that they obey the Cauchy-Riemann conditions for differentiability of a complex function. This in turn explains the conformal property of this mapping.

## 7.4 Relation to the Lemniscate Constant

The constant $K_e$ is actually related to another mathematical constant known as the *lemniscate constant*. The definition of the lemniscate constant [Langer 2011] is

$$L_e = \frac{\left[\Gamma\left(\frac{1}{4}\right)\right]^2}{2\sqrt{2\pi}}$$

where Γ is the gamma function. The numerical value of $L_e$ is approximately 2.62205755

The constant $K_e$ is related to $L_e$ by the following equation:

$$K_e \;=\; \frac{L_e}{\sqrt{2}} \;=\; \frac{\left[\Gamma\left(\frac{1}{4}\right)\right]^2}{4\sqrt{\pi}}$$

## 7.5 Verifying Conformality

One way to verify the conformality of the mapping is by looking at the mapping diagrams in page 6. The top figure shows a radial grid inside a circular disc mapped to a square. The radial lines and circular contours meet at 90° inside the circular disc. Since this mapping is conformal, the corresponding curves and squircle-like contours should also meet at 90° inside the square region. This can be verified visually by inspection.

Likewise, we can observe something similar for the bottom figure showing a rectangular grid inside a square region. This rectangular grid is composed of horizontal and vertical lines that meet at 90° inside the square. After mapping to the circular disc, the corresponding curves also meet at 90° inside the circle.

The conformal map between the circular disc and the square has many applications in science and engineering. For example, the Peirce quincuncial map projection [Fong 2011] used in geography and panoramic photography relies on this conformal map as an intermediate step.

## 7.6 Software Implementation

In order to implement this conformal mapping on a computer, one needs to be able to calculate special functions such as $F$ and $cn$. Source code for numerical computation of these special functions is readily available in open-source libraries such as Boost and the GNU Scientific Library. There is also a reference implementation appearing in the popular Numerical Recipes book [Press 1992]. There are well-established fast and robust algorithms for computing the $F$ and $cn$ special functions [Carlson 1977].

It is important to reiterate here that the Schwarz-Christoffel mapping requires the complex versions of the incomplete Legendre elliptic integral of the 1st kind $F$ and the Jacobi elliptic function $cn$. Both of these special functions are well defined for complex number arguments. L.M. Milne-Thompson provides the equations for the complex extension of these special functions in the classic AMS-55 reference book [Abramowitz 1964]

.

## 8. Derivation of the Fundamental Conformal Diagram

In this section, we will derive the fundamental conformal diagram shown in Section 7.1 from basic principles of the Schwarz-Christoffel transformation. First, let us review the equation for mapping the circular disc to a simple polygon. Driscoll and Trefethen's book [Driscoll 2002] is the definitive authority on fundamental concepts of Schwarz-Christoffel mappings. It provides the general formula for a conformal mapping between the circular disc and a simple polygon:

$$z \;=\; f(w) \;=\; \mathcal{A} + \mathcal{C} \int_0^w \prod_{j=1}^{n} \left(1 - \frac{\xi}{w_j}\right)^{\alpha_j - 1} \; d\xi$$

This equation is quite complicated and requires several definitions of the numerous variables in the equation. The complex variables $\boldsymbol{w}$ and $z$ represent points in the complex plane where the mapping occurs. Specifically, $\boldsymbol{w}$ represents points inside the unit circular disc centered at the origin. Meanwhile, $z$ represents points inside the target polygon. Here, the target polygon is n-sided with vertices located at $\boldsymbol{w}_j$ having vertex angles $\alpha_j$.

The scripted letters $\mathcal{A}$ and $\mathcal{C}$ are complex-valued constants that specify the position, orientation, and scale of the target polygon in the complex plane. We are free to choose any convenient value for these constants to control the placement and size of the target polygon. Note that $\mathcal{C}$ is a rotational and scaling constant with the restriction of being non-zero.

The mapping equation becomes much simpler if we restrict the target polygon to a regular polygon [Langer 2011], whereupon it simplifies to

$$z = f(w) = \; \mathcal{A} + \mathcal{C} \int_0^w \frac{d\tau}{(1 - \tau^n)^{2/n}}$$

Moreover, for our square case where *n=4*, the equation further simplifies to

$$z = f(w) = \; \mathcal{A} + \mathcal{C} \int_0^w \frac{d\tau}{\sqrt{(1 - \tau^4)}}$$

Now, we will segue our discussion to elliptic integrals. Specifically, we are interested in the incomplete Legendre integral of the 1st kind *F*. The Legendre normal form of this integral is given by NIST DLMF 19.2.4 as

$$F(\varphi, k) = \; \int_0^{\sin\varphi} \frac{d\tau}{\sqrt{1 - \tau^2}\sqrt{1 - k^2\tau^2}}$$

If we consider an imaginary modulus *k=i*, we can simplify to get

$$F(\varphi, i) = \int_0^{\sin\varphi} \frac{d\tau}{\sqrt{1 - \tau^2}\sqrt{1 - i^2\tau^2}} = \int_0^{\sin\varphi} \frac{d\tau}{\sqrt{1 - \tau^2}\sqrt{1 + \tau^2}} = \int_0^{\sin\varphi} \frac{d\tau}{\sqrt{(1 - \tau^4)}}$$

Using this, we can rewrite our Schwarz-Christoffel mapping equation in terms of a Legendre elliptic integral

$$z = f(w) = \mathcal{A} + \mathcal{C} \quad F(\sin^{-1} w\,, i)$$

The standard AMS-55 reference [Abramowitz 1964] provides an equation (17.4.17) to handle a Legendre elliptic integral with an imaginary modulus:

$$F(\xi, i) = \frac{1}{\sqrt{2}} F\left(\tfrac{\pi}{2}, \tfrac{1}{\sqrt{2}}\right) - \frac{1}{\sqrt{2}} F\left(\tfrac{\pi}{2} - \xi, \tfrac{1}{\sqrt{2}}\right)$$

Hence,

$$F(\sin^{-1} w\,, i) = \frac{1}{\sqrt{2}} F\left(\tfrac{\pi}{2}, \tfrac{1}{\sqrt{2}}\right) - \frac{1}{\sqrt{2}} F\left(\tfrac{\pi}{2} - \sin^{-1} w\,, \tfrac{1}{\sqrt{2}}\right)$$

Furthermore, using the trigonometric identity: $\cos^{-1} w = \frac{\pi}{2} - \sin^{-1} w$, this simplifies to

$$F(\sin^{-1} w, i) = \frac{1}{\sqrt{2}} F\left(\frac{\pi}{2}, \frac{1}{\sqrt{2}}\right) - \frac{1}{\sqrt{2}} F(\cos^{-1} w, \frac{1}{\sqrt{2}})$$

Recall the definition of the constant $K_e = F\left(\frac{\pi}{2}, \frac{1}{\sqrt{2}}\right)$. We can use this to simplify our mapping equation to

$$z = f(w) = \mathcal{A} + \mathcal{C}\left( \frac{1}{\sqrt{2}} K_e - \frac{1}{\sqrt{2}} F\left(\cos^{-1} w, \frac{1}{\sqrt{2}}\right) \right)$$

Here, we are free to choose whatever convenient values we want for the complex constants $\mathcal{A}$ and $\mathcal{C}$. For example, if we select $\mathcal{A} = 0$ and $\mathcal{C} = 1$, we get a square centered at the origin of the complex plane.

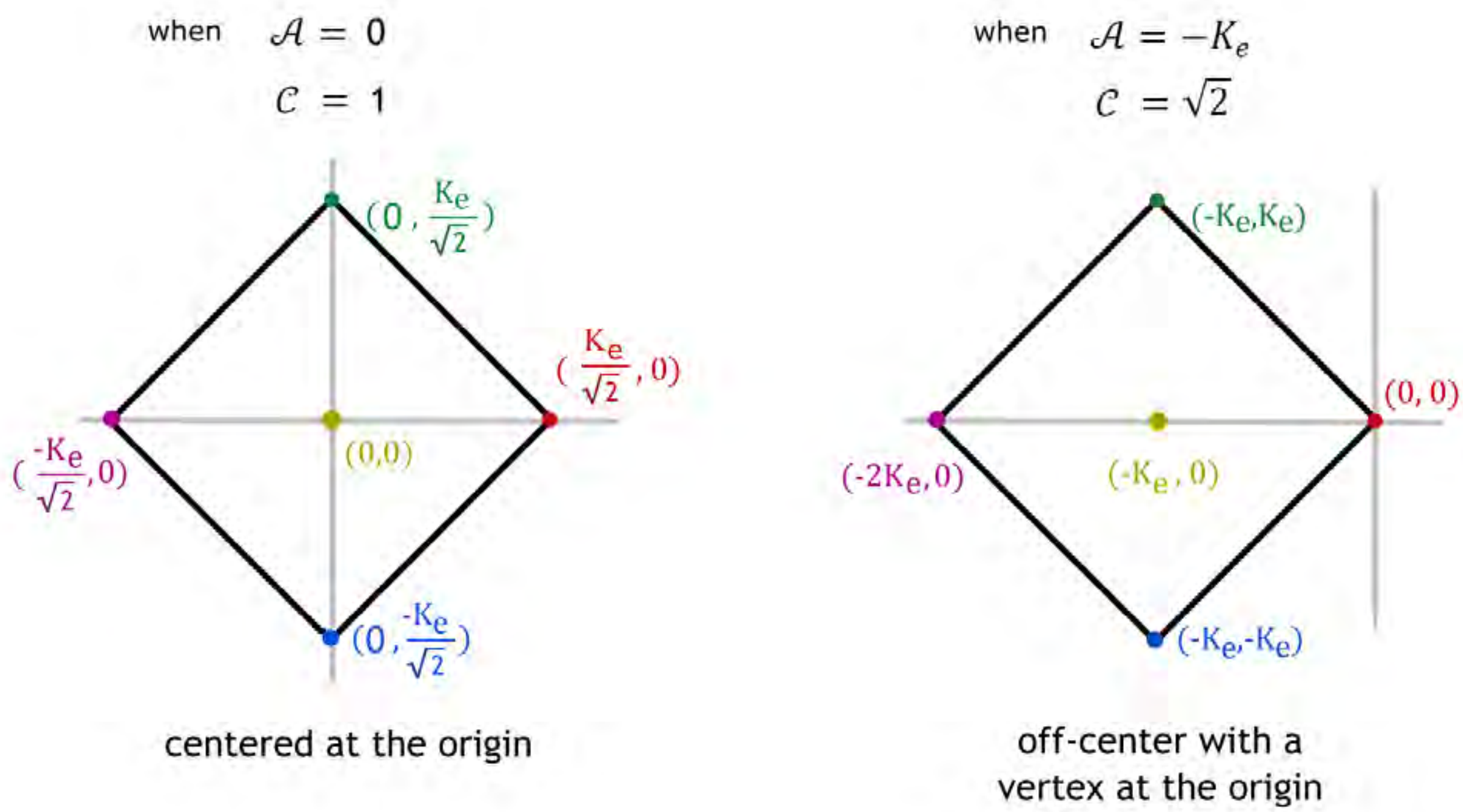

The four vertices of the square are scaled versions of the 4th roots of unity. The constant scale factor is $\frac{K_e}{\sqrt{2}}$. This arises from evaluating *f(w)* at *w=1*

$$f(1) = \frac{1}{\sqrt{2}} K_e - \frac{1}{\sqrt{2}} F\left(\cos^{-1} 1, \frac{1}{\sqrt{2}}\right) = \frac{1}{\sqrt{2}} K_e - \frac{1}{\sqrt{2}} F\left(0, \frac{1}{\sqrt{2}}\right) = \frac{K_e}{\sqrt{2}}$$

Meanwhile, it is convenient to select $\mathcal{A}$ and $\mathcal{C}$ so that we can get an even simpler mapping equation. Specifically, if we choose

$$\mathcal{A} = -K_e \qquad\qquad \mathcal{C} = \sqrt{2}$$

and plug these values into the mapping equation

$$z = f(w) = \mathcal{A} + \frac{\mathcal{C}}{\sqrt{2}} K_e - \frac{\mathcal{C}}{\sqrt{2}} F\left(\cos^{-1} w, \frac{1}{\sqrt{2}}\right)$$

we get

$$z = f(w) \; = \; -K_e + \; \frac{\sqrt{2}}{\sqrt{2}} K_e - \frac{\sqrt{2}}{\sqrt{2}} F\left(\cos^{-1} w \,, \tfrac{1}{\sqrt{2}}\right)$$

$$= \; -F\left(\cos^{-1} w \,, \tfrac{1}{\sqrt{2}}\right)$$

or equivalently,

$$-z = F\left(\cos^{-1} w \,, \tfrac{1}{\sqrt{2}}\right)$$

After inverting the elliptic integral, we get

$$cn\left(-z, \tfrac{1}{\sqrt{2}}\right) = w$$

Recall that the Jacobi elliptic function cn is an even function, so $cn(-z, k) = cn(z, k)$. This reduces our mapping equation to

$$\boldsymbol{w = cn\left(z, \tfrac{1}{\sqrt{2}}\right)}$$

Thus our fundamental conformal diagram arises from the choice of $\mathcal{A} = -K_e$ and $\mathcal{C} = \sqrt{2}$ in the Schwarz-Christoffel mapping equation. This completes the derivation ■.

As a side note, we would like to mention that the Schwarz-Christoffel mapping in our canonical mapping space partially arises from the choice of $\mathcal{A} = 1 - i$ and $\mathcal{C} = \frac{-\sqrt{-2i}}{K_e}$ in the mapping equation.

## 9. Conformal and Equiareal Non-Existence

In this section, we claim that there is no disc-to-square mapping in our canonical mapping space that is both conformal and equiareal. The argument is relatively simple and is based on the Riemann Mapping theorem from complex analysis.

### 9.1 Riemann Mapping Theorem

We briefly review the statement of the Riemann Mapping theorem here. The Riemann Mapping theorem states that if **U** is a non-empty simply connected open subset of the complex number plane $\mathbb{C}$, which is not all of $\mathbb{C}$, then there exists a bijective holomorphic mapping *f* from **U** onto the open unit disk. Moreover, the theorem states that the mapping is unique if we fix a point and the orientation of the mapping.

### 9.2 Uniqueness of Schwarz-Christoffel Mapping

For each of the disc to square mappings discussed in this paper, our canonical mapping space imposes some strict restrictions on the mappings. For example, we want the center of the circle to map to the center of the square, i.e. the *canonical central constraint*. Also, we have the *canonical axial constraints* from section 1.3 of this paper. These constraints essentially fix a point and the orientation of the mappings.

Using the Riemann Mapping theorem, we can claim that Schwarz-Christoffel mapping of the circle to the square is the only conformal mapping between these two shapes, subject to the restrictions of our canonical mapping space. In other words, our version of the Schwarz-Christoffel mapping is uniquely conformal.

### 9.3 The Schwarz-Christoffel Mapping is not Equiareal

In order to show the non-existence of a disc-to-square mapping that is both conformal and equiareal in our canonical mapping space, all we have to do is show that our version of the Schwarz-Christoffel mapping is not equiareal. Consider the Schwarz-Christoffel mapping diagram shown below.

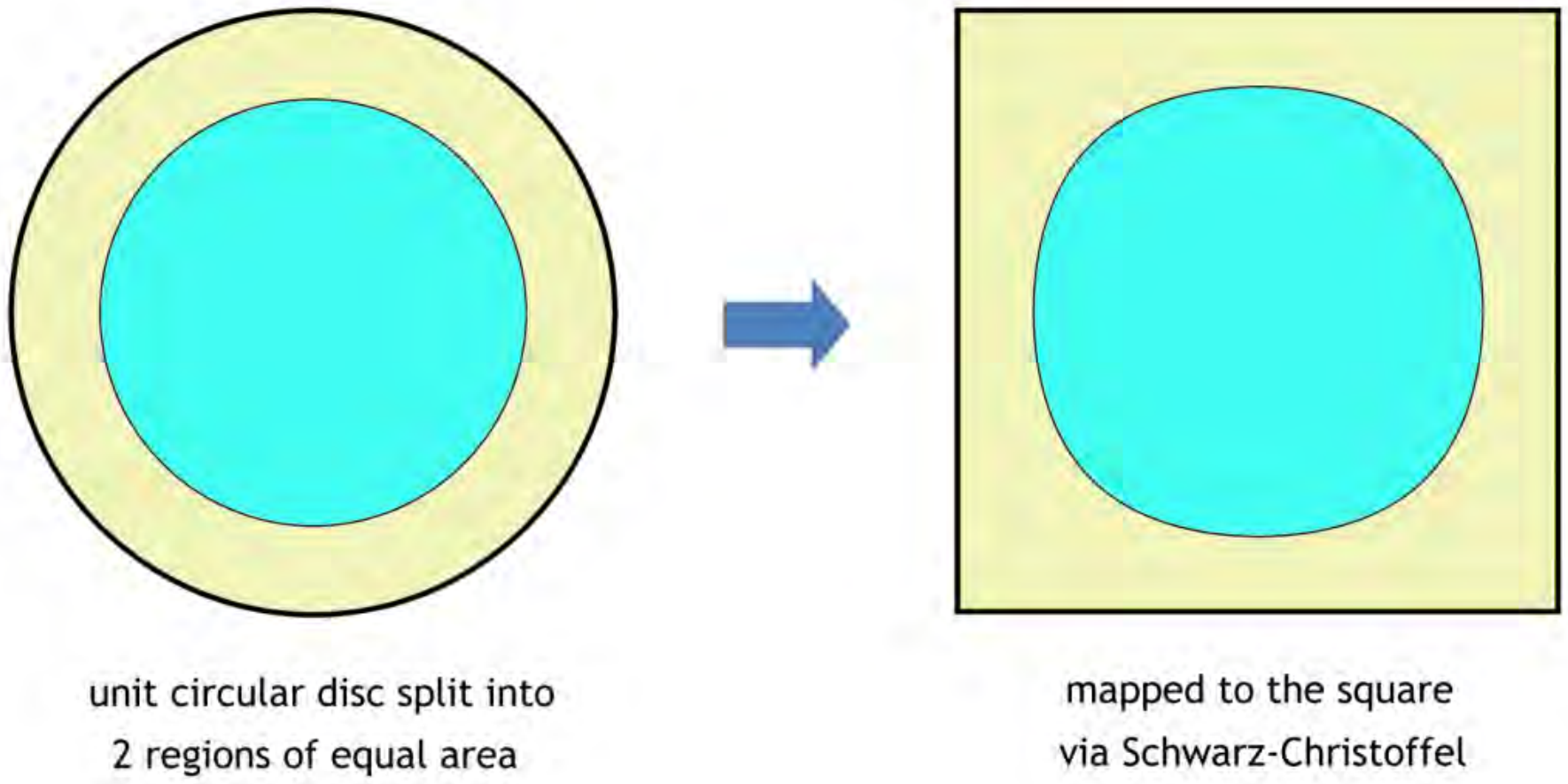

The unit circular disc on the left is split into two regions of equal area. These regions consist of a smaller central disc and an annular ring surrounding it. In order for both of these regions to have the same area, the smaller central disc must have a radius of $\frac{1}{\sqrt{2}}$. The inner circle enclosing the central disc can be expressed parametrically in the complex plane as $\frac{1}{\sqrt{2}}e^{i\tau}$ with $\tau \in [0, 2\pi]$.

The square on the right is the mapped version of the circular disc using the Schwarz-Christoffel mapping. This square has an area of 4. It is split into two regions – an inner region and an enclosing outer region. If the mapping is equiareal, then these two regions must have the same area. In other words, if $A_S$ is the area of the inner region, an equiareal mapping necessitates:

$$A_s = 4 - A_s = 2$$

We now attempt to calculate $A_S$ to check whether the mapping is equiareal. First, we have to find an equation for the enclosing curve of the inner region. Using the disc-to-square mapping equation in Section 7.3, this curve in the complex plane can written down parametrically as

$$\mathcal{S}(\tau) = 1 - i - \frac{\sqrt{-2i}}{K_e} F\left(\cos^{-1}\left(\frac{1}{\sqrt{2}} e^{i\tau}\sqrt{i}\right), \frac{1}{\sqrt{2}}\right)$$

for $\tau \in [0, 2\pi]$. Likewise, the *x* and *y* coordinates of this curve can be expressed parametrically as

$$x_s(\tau) = Re(1 - i - \frac{\sqrt{-2i}}{K_e} F\left(\cos^{-1}\left(\frac{1}{\sqrt{2}} e^{i\tau}\sqrt{i}\right), \frac{1}{\sqrt{2}}\right))$$

$$y_s(\tau) = Im(1 - i - \frac{\sqrt{-2i}}{K_e} F\left(\cos^{-1}\left(\frac{1}{\sqrt{2}} e^{i\tau}\sqrt{i}\right), \frac{1}{\sqrt{2}}\right))$$

Using these parametric equations, we can calculate the area $A_S$ of the inner region. First recall the formula for area in polar coordinates:

$$A = \frac{1}{2}\int_a^b r^2 \, d\theta$$

In polar coordinates, we have

$$r_s(\tau) = \sqrt{\left(x_s(\tau)\right)^2 + \left(y_s(\tau)\right)^2}$$

so the equation for the area of the inner region is

$$A_s = \frac{1}{2}\int_0^{2\pi} [\, \left(x_s(\tau)\right)^2 + \left(y_s(\tau)\right)^2 \,] \, d\tau$$

This definite integral is quite complicated but it can be calculated numerically using standard mathematical software packages or through custom-written numerical analysis software. In fact, we did this and found

$$A_s \approx 1.83 \quad \neq 4 - A_s$$

Hence, our Schwarz-Christoffel mapping cannot be equiareal. This shows that subject to our canonical mapping space restrictions, there is no disc-to-square mapping that is both conformal and equiareal.

# 10. Applications

## 10.1 The Conformal Square and Hyperbolic Art

The Poincare disk is one of the most interesting models of non-Euclidean hyperbolic geometry. In fact, this model of the hyperbolic plane has inspired artwork such as M.C. Escher's circle limit woodcuts. Using the different mappings discussed in this paper, one can convert the Poincare disc to a square. The Schwarz-Christoffel mapping is the most appropriate for this task because of its conformal nature.

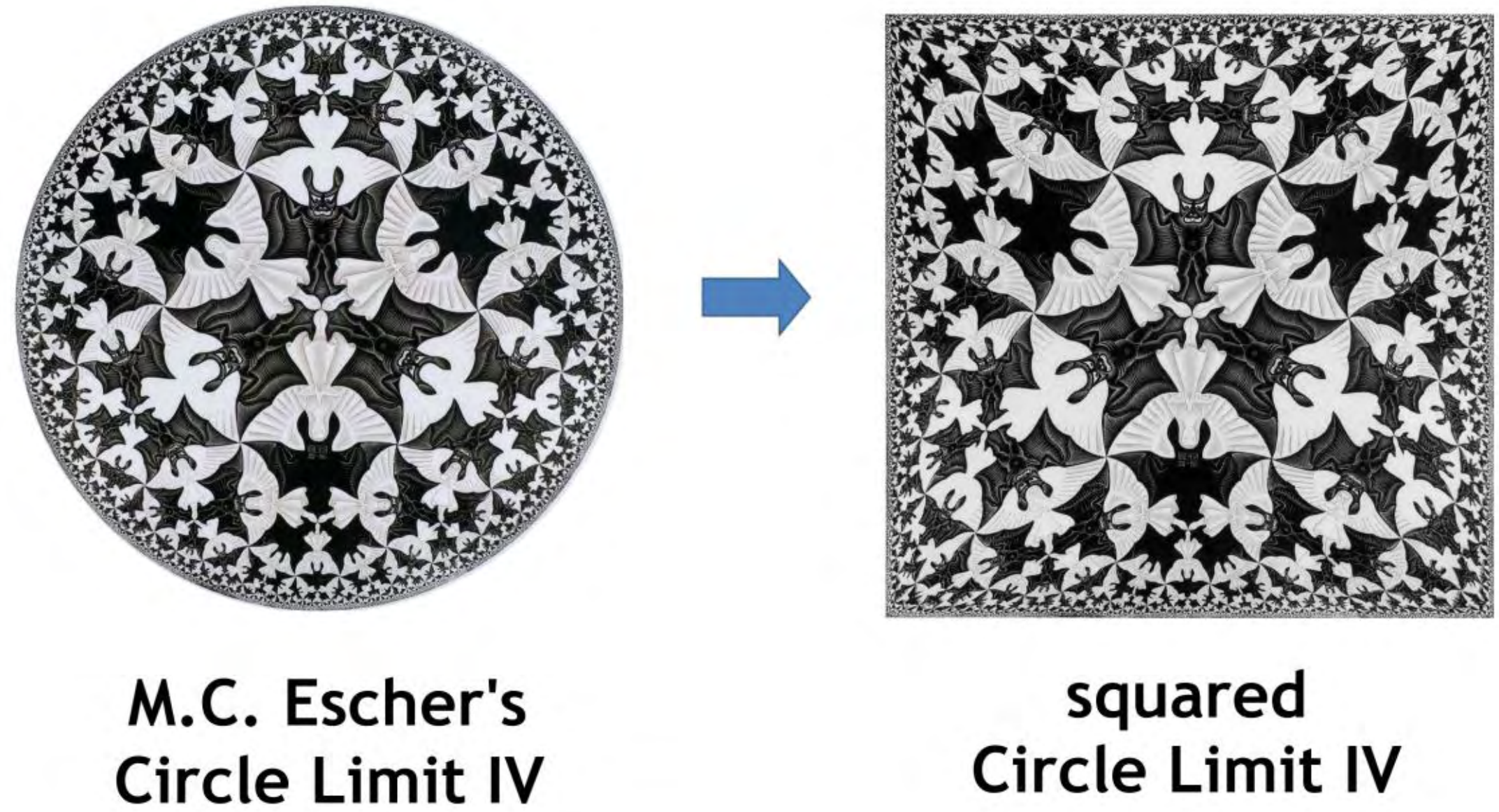

Herbert Müller [Müller 2013] has used Schwarz-Christoffel transformations on the complex plane to convert M.C. Escher's Circle Limit IV to polygons with different number (n) of sides. His paper shows results for $n=3$ and $n=6$ but not for $n=4$ (the square).

## 10.2 Logo Design and Artwork

The circle and the square are very common shapes used in logos. It is certainly useful to have methods to convert from one shape to the other as part of the designer's toolbox. To illustrate this, we convert the Fields medal, which is circular in shape, to a square one below. Of course, the Fields medal awards ceremony is a highlight of the quadrennial International Congress of Mathematicians (ICM) for which this work was shown in 2014.

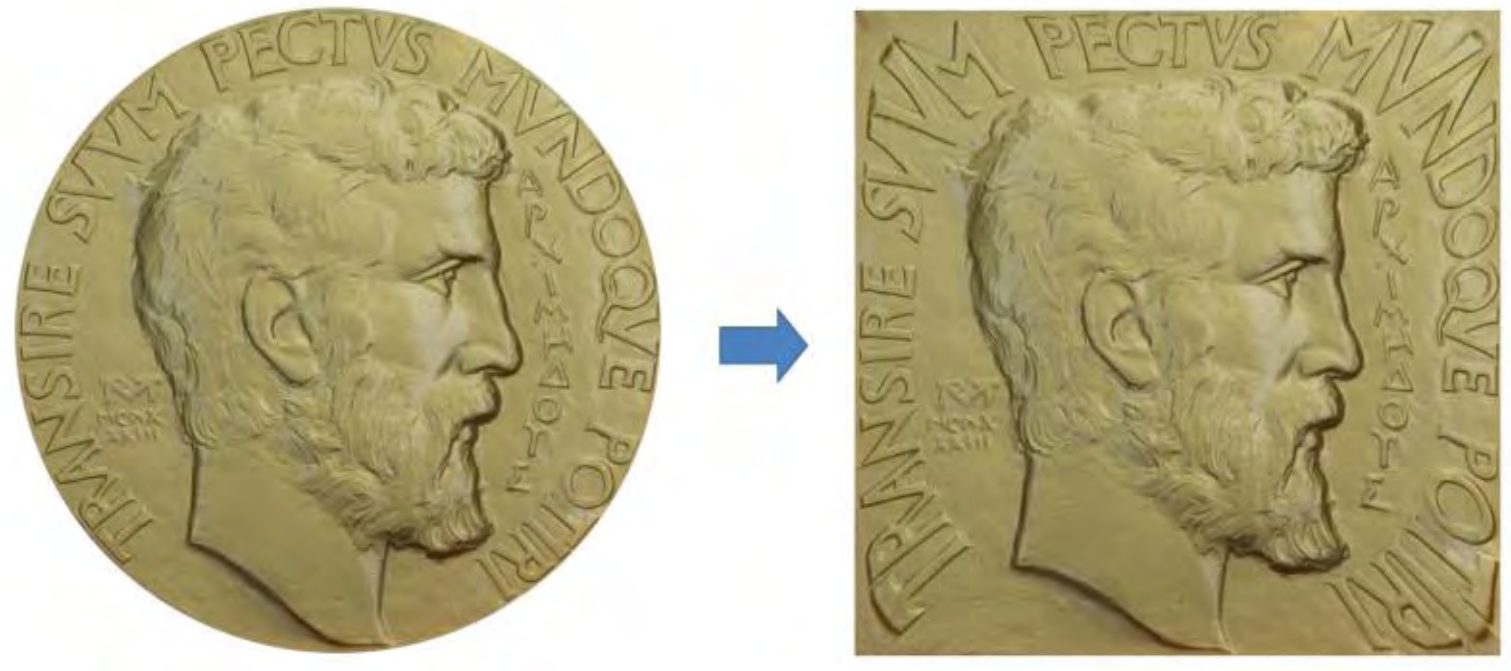

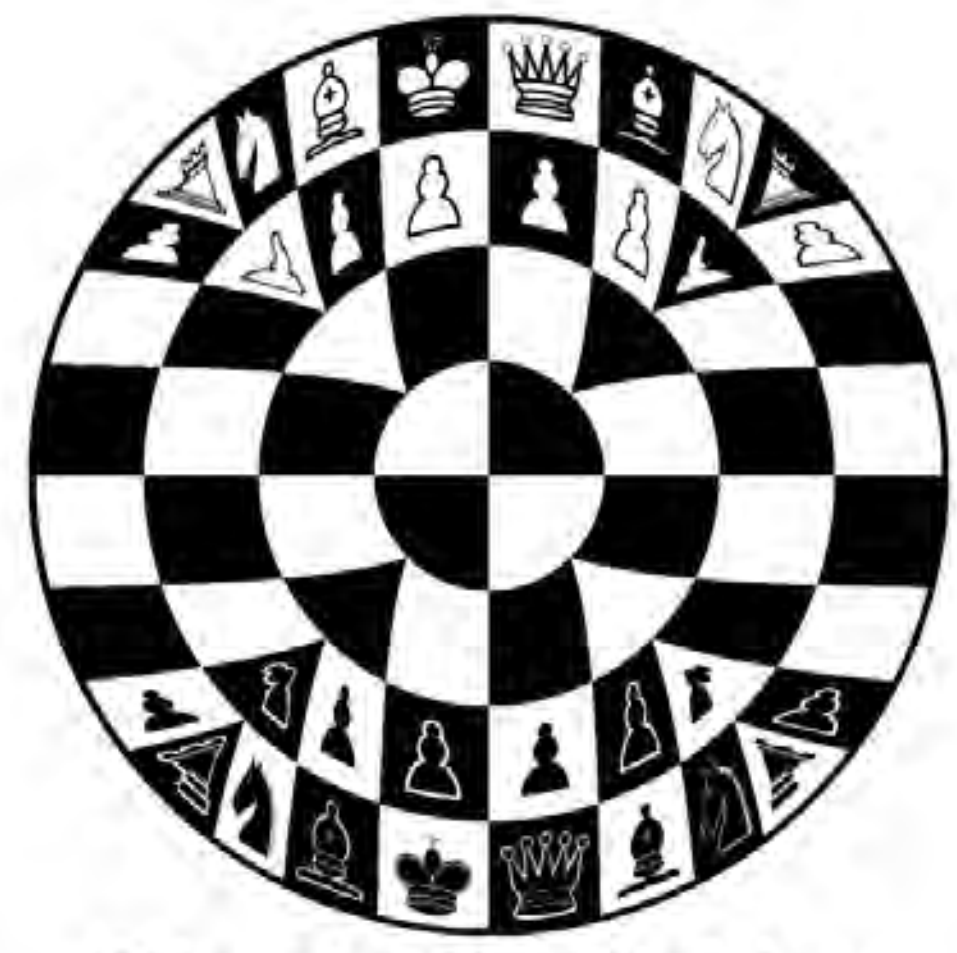

Simple Stretching

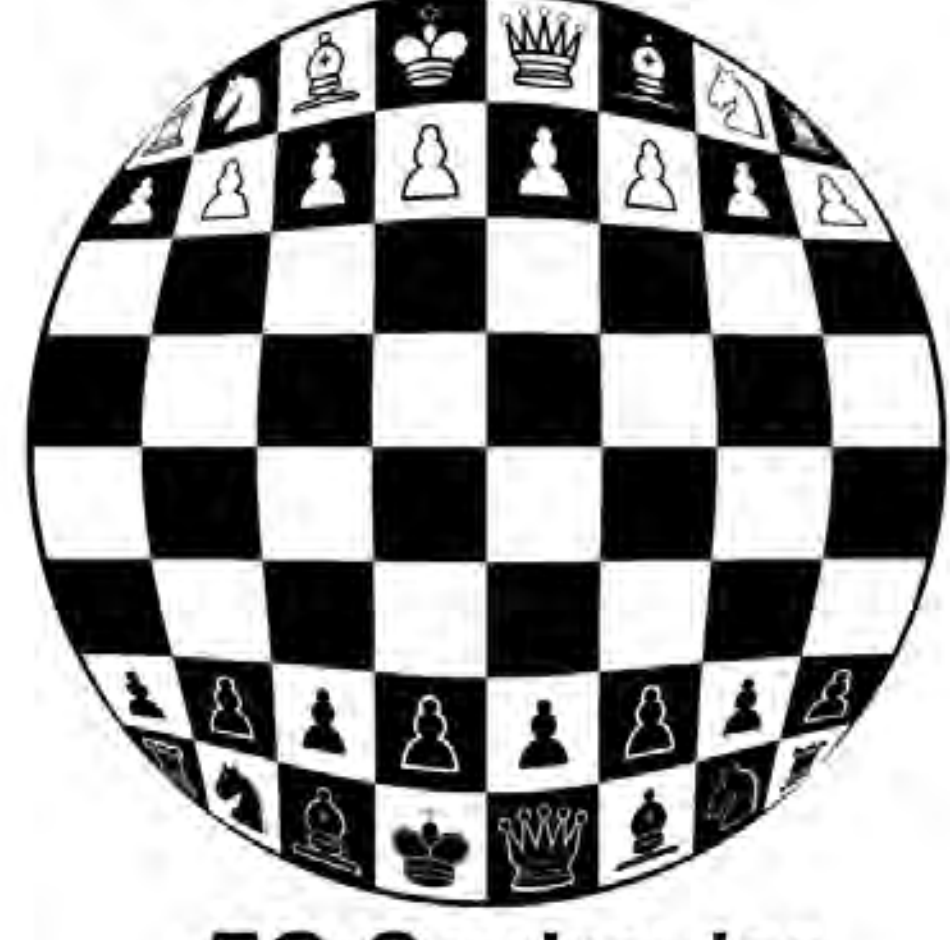

FG-Squircular

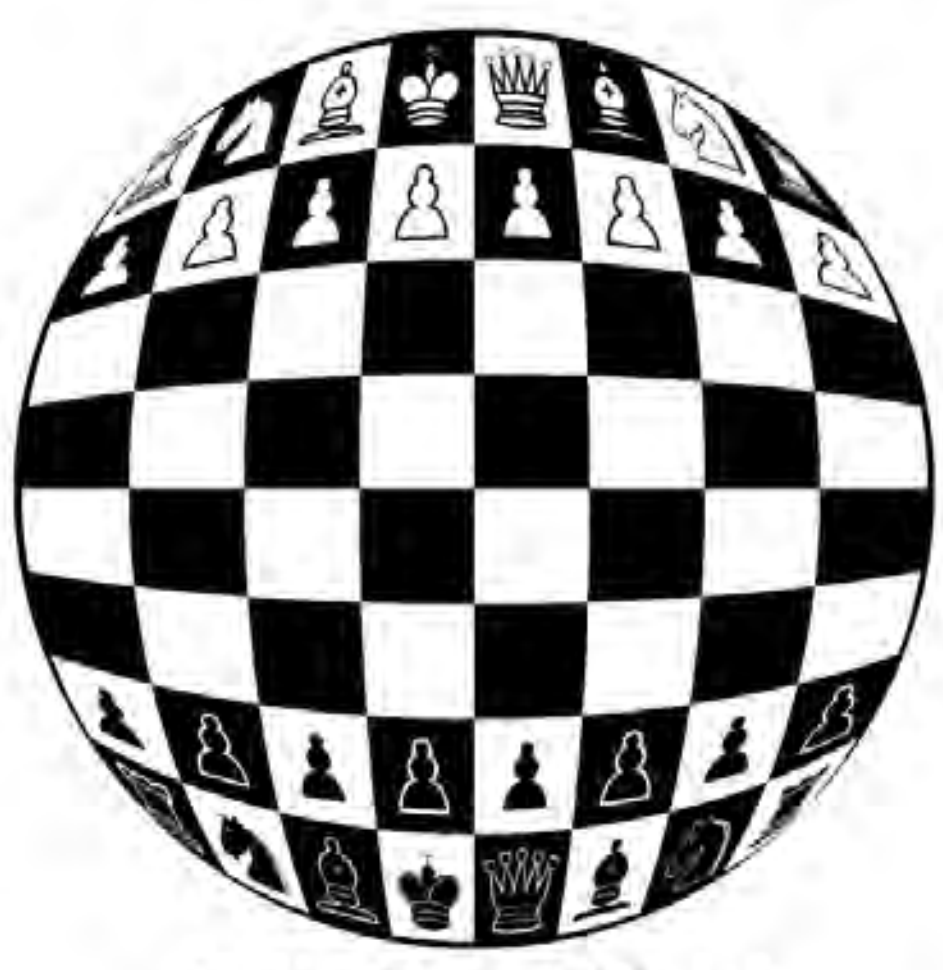

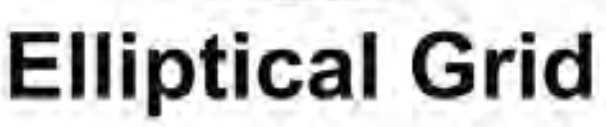

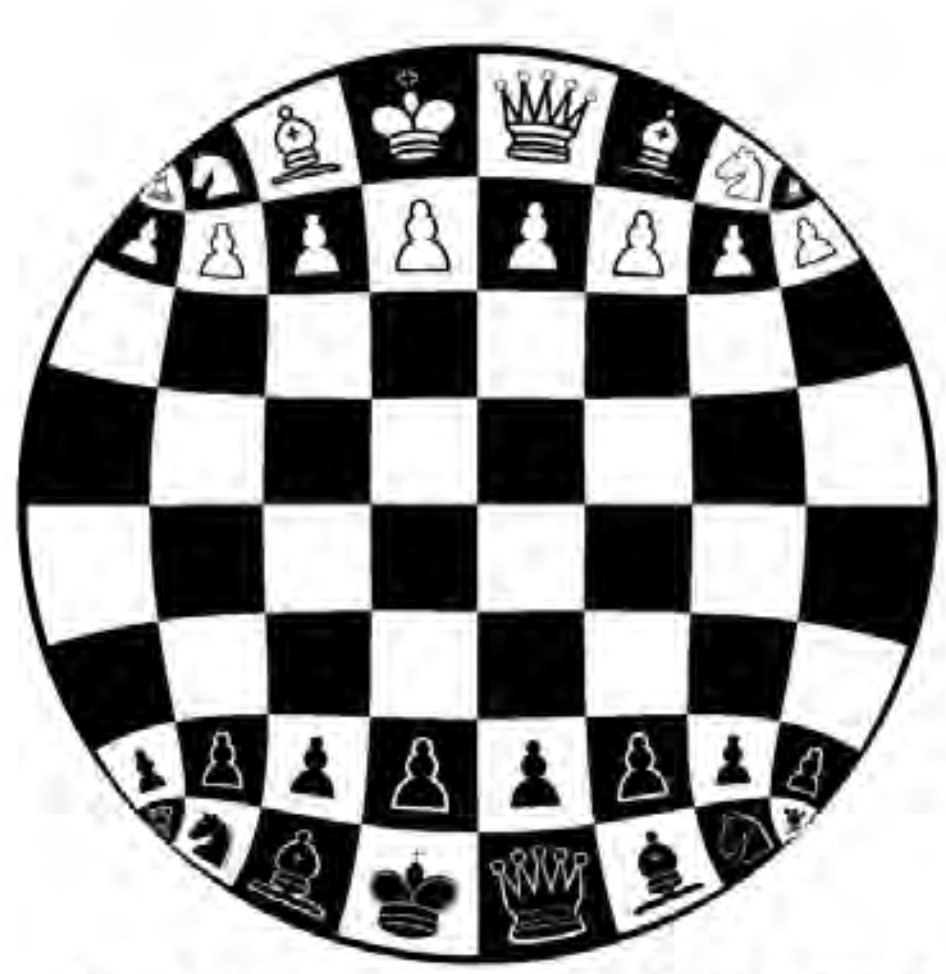

Schwarz-Christoffel

The four different mappings presented in this paper can be used to convert square-shaped designs into circular artwork. For example, the figure above show results from mapping a square chessboard into a circular disc. Each of the four mappings produces quite different results and distortion characteristics. These distortions are most apparent at the four corners of the board where the rooks are located.

It is quite evident that the Simple Stretching map has unsightly artifacts along its diagonals. In fact, we only included it in this paper because of its simplicity. It is largely unsuitable for our intended applications.

The Schwarz-Christoffel mapping is conformal. This is quite noticeable on its corresponding circular map. By and large, each of the chess pieces on the board have no shape distortion and very much resemble the shapes on the square chessboard. However, there is significant size distortion in the mapping. This distortion is most prominent with the corner rooks. The rooks are considerably smaller than the other chess pieces.

The Elliptical Grid mapping produces a nice and uniform checkered pattern on the circular disc. Horizontal and vertical lines on the square chessboard map into a curvilinear grid consisting of elliptical arcs on the circular disc. This characteristic definitely makes the Elliptical Grid mapping quite desirable. However, there is considerable shape distortion at the four corners of the chessboard. In fact, it is quite noticeable that the rooks are horribly stretched and bent out of shape.

The FG-Squircular mapping offers a good compromise between shape and size distortion. The four corner rooks are still noticeably smaller than the other chess pieces, but their shapes are not as deformed as in the Elliptical Grid mapping. The rooks also appear much larger than their counterparts in the Schwarz-Christoffel mapping.

## 10.3 Azimuthal Panoramas

Azimuthal projections of spherical panoramas produce naturally circular images. However, since most of the world's photographs are rectangular in shape, it is desirable to convert these circular images to square ones. An example azimuthal panorama is provided below along with its mapping to a square [Fong 2014]. The FG-Squircular mapping was used for the mapping because of its radial nature. The radial property makes it suitable for extending the radial grid of the azimuthal projection to the square.

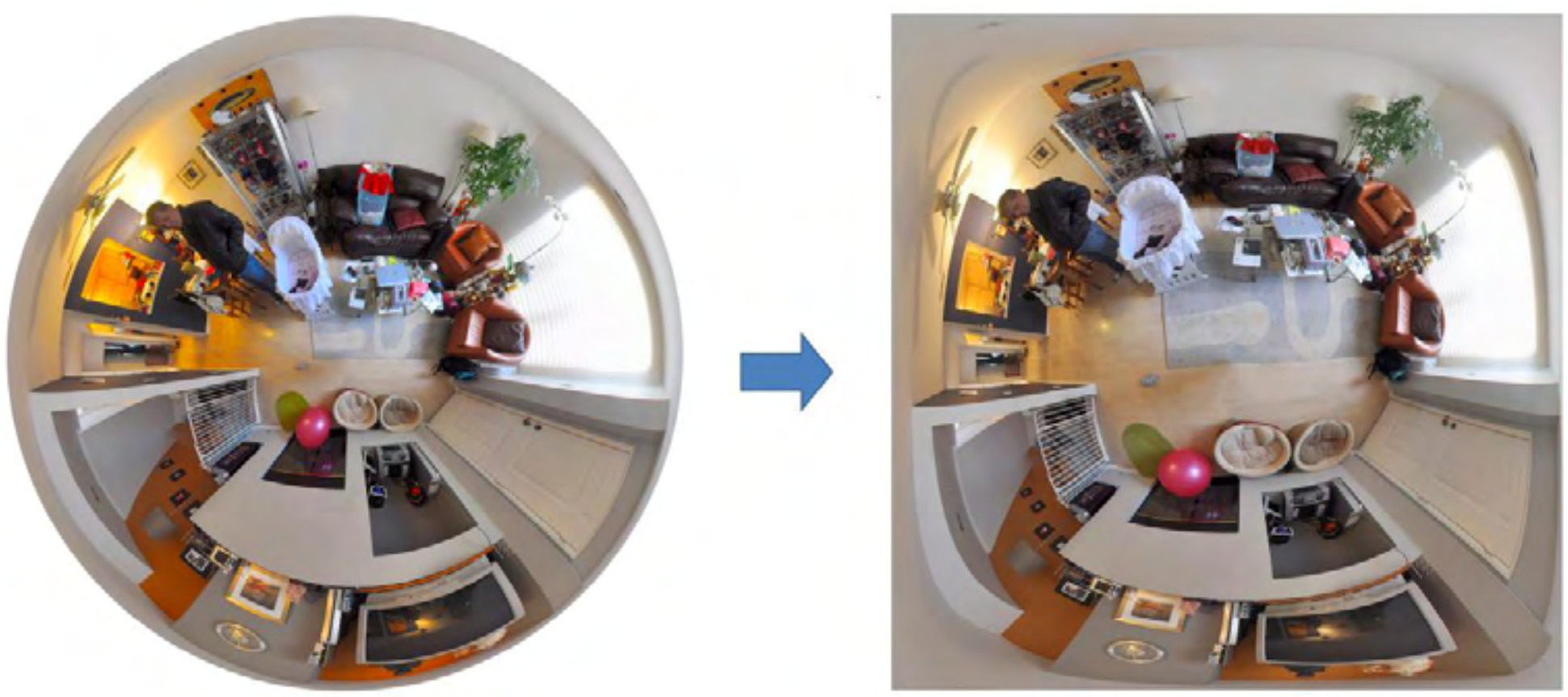

## 10.4 Quick and Dirty Defishing

There are two types of fisheye lens – frame-filling rectangular fisheye lens and circular fisheye lens. Photographs taken using circular fisheye lens are circular in shape. Photographers usually want to convert these circular photographs to square ones in order to straighten out curved lines distorted by the fisheye lens. This process is known as defishing.

Fisheye lens produce distortions known in the optics literature as barrel distortion [Falk 1986]. Visually speaking, this distortion magnifies the center of the image and falls off as you move away from the center. Qualitatively, this distortion looks quite similar to what the Elliptical Grid mapping does to a rectangular grid as shown in page 5. Consequently, it seems viable to use the Elliptical Grid equations to reverse this distortion.

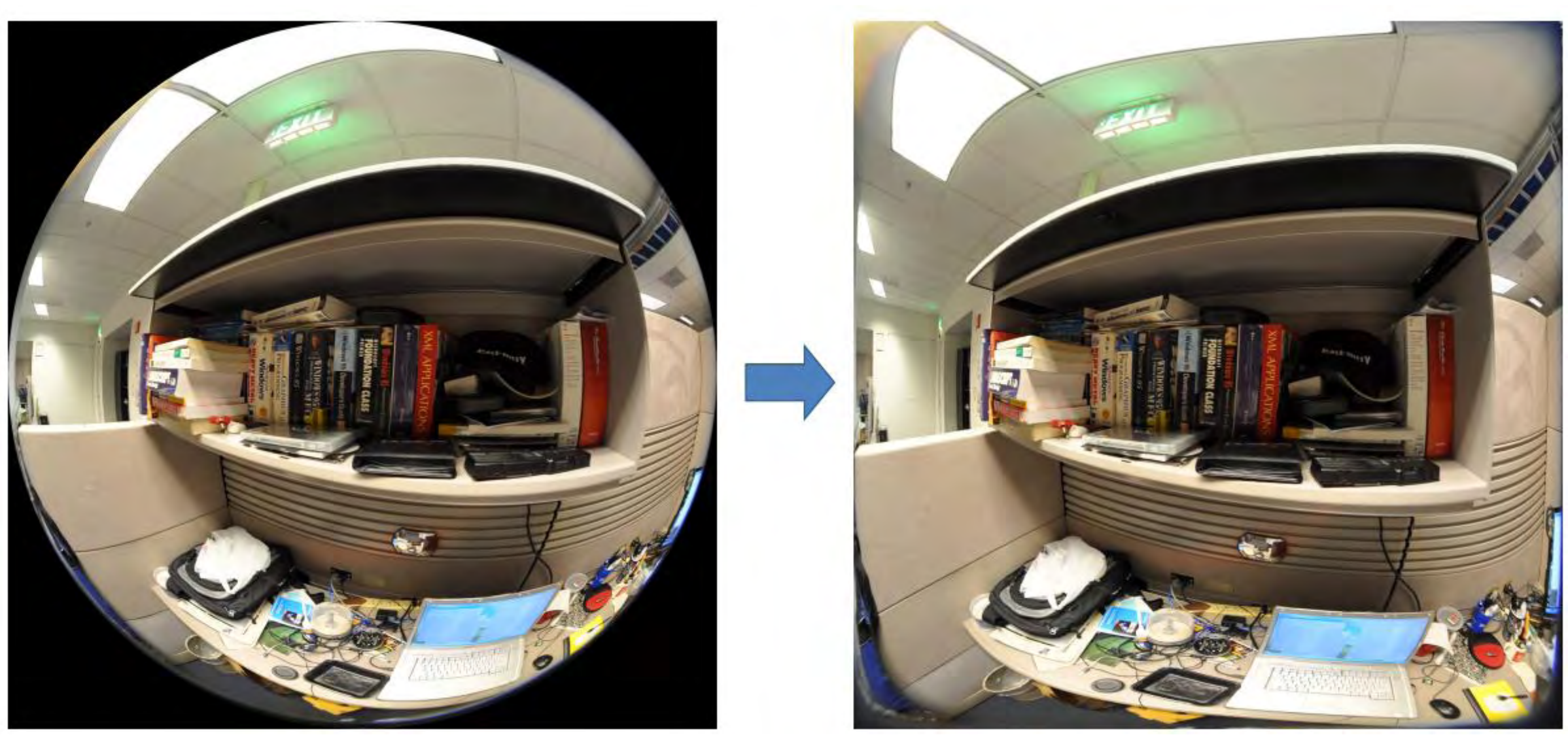

Using the disc to square mappings discussed in this paper, one can do quick and dirty defishing of photographs. We call it "quick and dirty" because the real defishing process is dependent on the optics of the fisheye lens in order to be done correctly. There are actual physical models and equations for barrel distortion in lens [Bailey 2002].The mappings provided in this paper have no knowledge of the inner workings of fisheye lens optics, and hence give only a poor man's approximation of the defishing process.

Empirically, we have found that the Elliptical Grid mapping works well for defishing when compared to the other mappings discussed with this paper. This is probably because of the qualitative similarity between its distortion and the barrel distortion produced by fisheye lens. The sample image shown above has before and after photographs of the defishing process using the Elliptical Grid mapping.

## 11. Summary

In this paper, we presented and discussed four different mappings for converting a circular disc to a square region and vice versa. Each of the mappings has different properties and characteristics. The table below summarizes these properties.

| mapping | key property | contours inside square | notes/comments |
|---|---|---|---|
| Simple Stretching | radial | squares | Circular disc is just linearly stretched to a square |
| FG-Squircular | radial | FG-squircles | Good compromise between equiareal and conformal |
| Elliptical Grid | grid of elliptical arcs | FG-squircles | Rectangular grid mapped to elliptical grid inside circle |
| Schwarz-Christoffel | conformal | squircle-like curve | Considerable size distortions at the four corners |

We also discussed different applications of the mappings between the disc and the square. We thereupon showed pictures of the mappings used for logo design, panoramic photography, and art. In addition, we showed that the elliptical grid mapping can be used for defishing photographs with barrel distortion from fisheye lens.

## 12. Acknowledgements

The author would like to thank Anton Shterenlikht, Ravi Dattatreya, Brian Vogel, Bill Reinke and other anonymous reviewers for providing insightful comments on this paper. The author would also like to give special thanks to Elaine Cheong for her invaluable conference travel tips during the presentation of the results in this paper.

# Addendum (2019)

*Even More Mappings of the Circular Disc to a Square*

In addition to the four discussed in this paper, the author has derived more mappings for squaring the circular disc. These additional mappings are discussed in the follow-up paper "*Elliptification of Rectangular Imagery*", which is available in

https://arxiv.org/abs/1709/07875

For example, the follow-up paper includes a variant of the FG-Squircular mapping called the 2-Squircular mapping which has much simpler equations.

## 2-Squircular Mapping

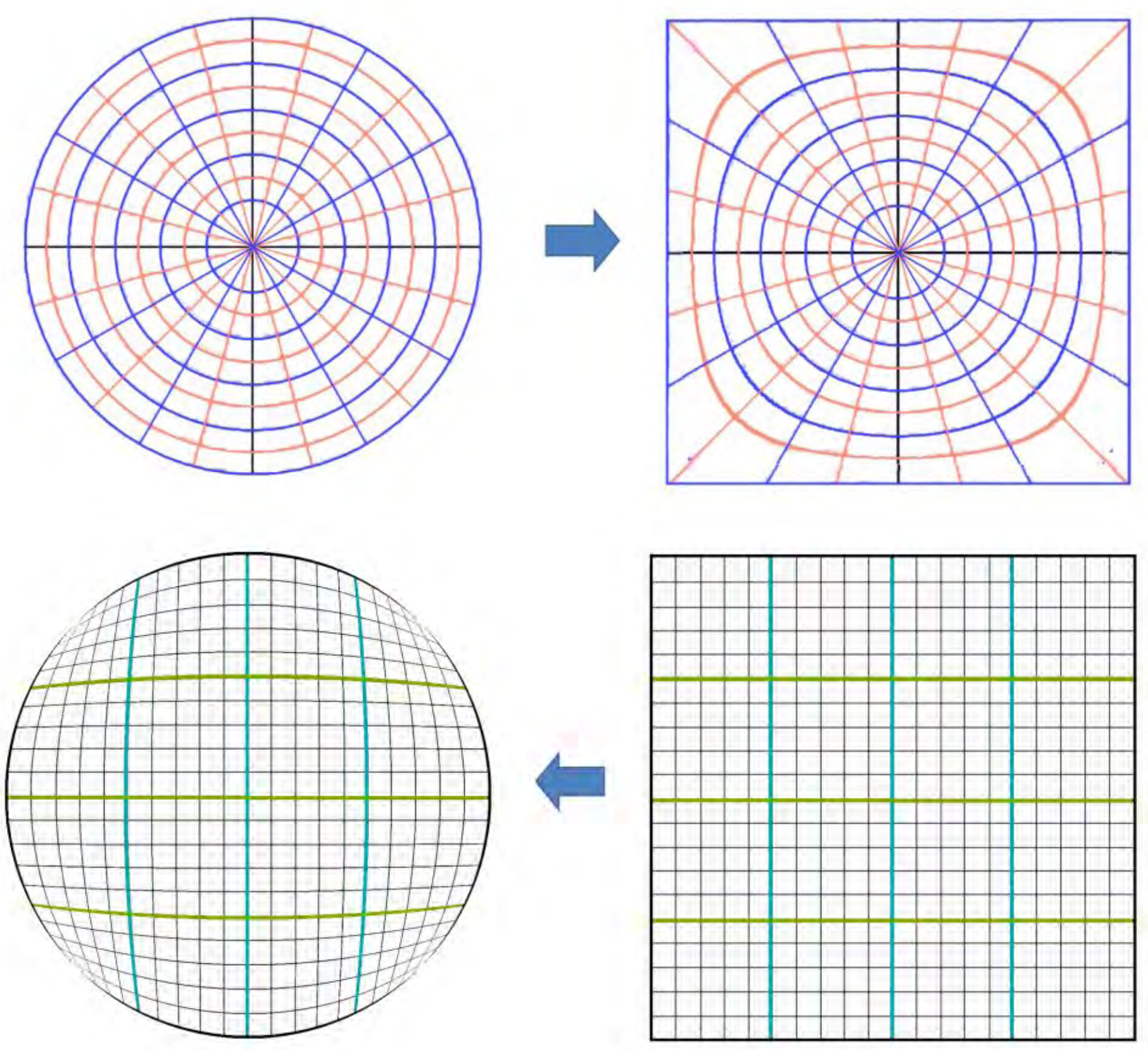

Disc to square mapping:

$$x = \frac{sgn(uv)}{v\sqrt{2}}\sqrt{1-\sqrt{1-4u^2v^2}}$$

$$y = \frac{sgn(uv)}{u\sqrt{2}}\sqrt{1-\sqrt{1-4u^2v^2}}$$

Square to disc mapping:

$$u = \frac{x}{\sqrt{1+x^2y^2}} \qquad\qquad v = \frac{y}{\sqrt{1+x^2y^2}}$$